\documentclass[11pt]{amsart}
\usepackage{amsmath}
\usepackage{amsxtra}
\usepackage{amscd}
\usepackage{amsthm}
\usepackage{amsfonts}
\usepackage{amssymb}
\usepackage{eucal}

\newtheorem{prop}{Proposition}

\newtheorem{claimm}{Claim}

\newtheorem{thmm}{Theorem}

\theoremstyle{remark}

\numberwithin{equation}{section}

\numberwithin{cor}{subsection}
\numberwithin{propconstr}{subsection}
\numberwithin{prop}{subsection}
\numberwithin{lem}{subsection}

\newcommand{\thmref}[1]{Theorem~\ref{#1}}
\newcommand{\secref}[1]{Sect.~\ref{#1}}
\newcommand{\lemref}[1]{Lemma~\ref{#1}}
\newcommand{\thmmref}[1]{Theorem~\ref{#1}}
\newcommand{\propref}[1]{Proposition~\ref{#1}}

\newcommand{\claimmref}[1]{Claim~\ref{#1}}

\newcommand{\corref}[1]{Corollary~\ref{#1}}
\newcommand{\conjref}[1]{Conjecture~\ref{#1}}
\newcommand{\remref}[1]{Remark~\ref{#1}}

\newcommand{\nc}{\newcommand}
\nc{\ssec}{\subsection}
\nc{\sssec}{\subsubsection}
\nc{\on}{\operatorname}

\nc{\Z}{{\mathbb Z}}
\nc{\NN}{{\mathbb N}}
\nc{\FF}{{\mathbf F}}
\nc{\CC}{{\mathbb C}}
\nc{\GG}{{\mathbb G}}
\nc\one{{\mathbf 1}}
\nc{\OO}{{\mathcal O}}
\nc{\DD}{{\mathbb D}}
\renewcommand{\AA}{{\mathbb A}}
\nc{\Fq}{{\mathbb F}_q}
\nc{\Fqb}{\overline{\mathbb F}_q}
\nc{\Ql}{\ol{\mathbb Q}_\ell}

\nc{\Eis}{\on{Eis}}
\nc{\Aut}{\on{Aut}}
\nc{\Rep}{\on{Rep}}
\nc{\Hom}{\on{Hom}}
\nc{\Loc}{\on{Loc}}
\nc{\Bun}{\on{Bun}}
\nc{\IC}{\on{IC}}
\nc{\Can}{\on{Can}}
\nc{\rk}{\on{rk}}
\nc{\Sh}{\on{Sh}}
\nc{\Perv}{\on{Perv}}
\nc{\Conv}{\on{Conv}}
\nc{\BunBb}{\overline{\Bun}_B}
\nc{\BunNft}{\Bun_N^{\F_T}}
\nc{\BunNftb}{\overline{\Bun}_N^{\F_T}}
\nc{\BunNftw}{\widetilde{\Bun}_N^{\F_T}}
\nc{\BunPb}{\overline{\Bun}_P}
\nc{\BunPbw}{\widetilde{\Bun}_P}
\nc{\GUb}{\overline{G/U}}
\nc{\Psib}{\overline{\Psi}}
\nc{\Psio}{\overset{o}{\Psi}}
\nc{\GUPb}{\overline{G/U(P)}}
\nc{\GPPb}{\overline{G/[P,P]}}
\nc\hl{h^{\leftarrow}}
\nc\hr{h^{\rightarrow}}
\nc\Hl{H^{\leftarrow}}
\nc\Hr{H^{\rightarrow}}
\renewcommand\i{{\mathfrak i}}
\nc\I{{\mathcal I}}

\nc{\F}{{\mathcal F}}
\nc{\A}{{\mathcal A}}
\nc{\W}{{\mathcal W}}
\nc{\G}{{\mathcal G}}
\nc{\J}{{\mathcal J}}
\nc{\Wb}{\overline{\W}}
\nc{\M}{{\mathcal M}}
\nc{\N}{{\mathcal N}}
\nc{\Y}{{\mathcal Y}}
\nc{\E}{{\mathcal E}}
\nc{\D}{{\mathcal D}}
\nc\Dh{\widehat{\D}}
\renewcommand{\O}{\hat{\mathcal O}}
\nc{\C}{{\mathcal C}}
\nc{\K}{\hat{\mathcal K}}
\renewcommand{\H}{{\mathcal H}}
\renewcommand{\S}{{\mathcal S}}
\nc{\T}{{\mathcal T}}
\renewcommand{\L}{{\mathcal L}}
\nc{\Gr}{\on{Gr}}

\nc{\p}{\overline{\mathfrak p}}
\nc{\q}{\overline{\mathfrak q}}
\nc{\qo}{{\mathfrak q}}
\nc{\po}{{\mathfrak p}}
\nc{\s}{{\mathfrak s}}

\nc{\pw}{\widetilde{\mathfrak p}}
\nc{\qw}{\widetilde{\mathfrak q}}
\nc{\jw}{\widetilde j}

\nc{\grb}{\overline{\Gr}}

\nc{\lambdach}{\check\lambda}
\nc{\Lambdach}{\check\Lambda}
\nc{\much}{\check\mu}
\nc{\omegach}{\check\omega}
\nc{\nuch}{\check\nu}
\nc{\etach}{\check\eta}
\nc{\alphach}{\check\alpha}
\nc{\betach}{\check\beta}
\nc{\rhoch}{\check\rho}
\nc\xl{\overline{x}}
\nc\yl{\overline{y}}
\nc\nul{\overline{\nu}}
\nc\mul{\overline{\mu}}
\nc\lambdal{\overline{\lambda}}
\nc\zerol{\overline{0}}

\nc{\Hb}{\overline{\H}}

\nc{\mc}{\mathcal}
\nc{\ga}{\gamma}
\nc{\GL}{{}^L G}
\nc{\PL}{\Lambda}
\nc{\la}{\lambda}
\nc{\arr}{\rightarrow}
\nc{\ol}{\overline}
\nc{\al}{\alpha}
\nc{\De}{\Delta}
\nc{\Gaf}{{\mathbb G}_{a,\Fq}}
\nc{\wt}{\widetilde}
\nc{\bi}{\bibitem}
\nc{\HH}{{\mathbf H}}
\nc{\WW}{{\mathbf W}}
\nc{\bn}{\BunNft}
\nc{\map}{\kappa}
\nc{\wh}{\widehat}
\nc{\bnx}{_{x,\infty}\BunNftb}
\nc{\canG}{\mathcal G}
\nc{\canB}{\mathcal B}
\nc{\canN}{{\mathcal N}^\epsilon}
\nc{\CanN}{\wt{\mathcal N}^\epsilon}
\nc{\canNnu}{{\mathcal N}^{\epsilon_\nu}}
\nc{\canNnup}{{\mathcal N}^{\epsilon_{\nu'}}}
\nc{\CanNnu}{\wt{\mathcal N}^{\epsilon_\nu}}
\nc{\CanNnup}{\wt{\mathcal N}^{\epsilon_{\nu'}}}
\nc{\out}{\on{out}}
\nc{\ft}{{\mathcal F}_T}
\nc{\zl}{\overline{z}}
\nc{\lal}{\overline{\la}}
\nc{\BunNftbm}{\ol{\Bun}_{N,\mu}^{\ft}}
\nc{\BunNftm}{\wt{\Bun}^{\F_T}_{N,\mu}}
\nc{\rhoc}{\check{\rho}}
\nc{\tboxtimes}{\widetilde{\boxtimes}}

\nc{\V}{{\mc V}}
\nc{\canNK}{\ol{N(K_x)}^{\epsilon_\nu}}
\nc{\jj}{{\mathfrak j}}
\renewcommand{\setminus}{-}

\title[Whittaker patterns on moduli spaces of
bundles]{Whittaker patterns in the geometry of moduli spaces of
bundles on curves}

\author{E. Frenkel}

\address{Department of Mathematics, University of California,
Berkeley, CA 94720, USA}

\author{D. Gaitsgory}

\address{Department of Mathematics, Harvard University, Cambridge, MA
02138 and School of Mathematics, Institute for Advanced Study,
Princeton, NJ 08540, USA}

\author{K. Vilonen}

\address{Department of Mathematics, Northwestern University, Evanston,
IL 60208, USA}

\date{July 1999; Revised: Aug 1999, Nov 2000}

\begin{document}

\maketitle


\section*{Introduction}

Whittaker functions are ubiquitous in the modern theory of automorphic
functions. In this paper we develop a geometric version of the local
Whittaker theory in the non--archimedian case. By this we mean finding
appropriate Whittaker sheaves, such that the functions obtained by
taking the traces of the Frobenius on the stalks are the local
Whittaker functions. The desire to construct such sheaves is motivated
by the geometric Langlands correspondence (see
\cite{Dr,La1,La2,previous}). Since Whittaker functions never have
compact support, we cannot obtain them from sheaves on
finite-dimensional varieties in an obvious way, and this makes the
problem of constructing Whittaker sheaves non-trivial. Our approach is
to use global methods, namely, the Drinfeld compactification of the
moduli stack of $N$--bundles.

More concretely, let $G$ be a (split connected) reductive group over
$\Fq$, and $N$ be its maximal unipotent subgroup. Consider the group
$G(\K)$ over the local field $\K=\Fq((t))$, and its maximal compact
subgroup $G(\O)$, where $\O = \Fq[[t]]$. The quotient $\Gr =
G(\K)/G(\O)$ can be given a structure of an ind-scheme over $\Fq$,
which is called the affine Grassmannian of $G$. Unramified Whittaker
functions give rise to functions on $\Gr$, which are
$N(\K)$--equivariant against a non-degenerate character
$\chi$. According to the general philosophy of
``faisceaux--fonctions'' correspondence, these functions should have
geometric counterparts as perverse sheaves on $\Gr$. The geometric
counterparts of the Whittaker functions should be sheaves on $\Gr$,
which are $N(\K)$--equivariant against $\chi$. It is natural to expect
that the appropriately defined category of such sheaves should reflect
the properties of the Whittaker functions. Unfortunately, the
$N(\K)$--orbits in $\Gr$ are infinite-dimensional and so there is no
obvious way to define this category on $\Gr$.

On the other hand, let $\Bun_N$ be the algebraic stack classifying
$N$--bundles on a smooth projective curve $X$ over $\Fq$.  V.~Drinfeld
has introduced a remarkable partial compactification $\ol{\Bun}_N$ of
$\Bun_N$.  In this paper we argue that an appropriate ``Whittaker
category'' can (and perhaps should) be defined on $\ol{\Bun}_N$
instead of $\Gr$. In fact, $\ol{\Bun}_N$ can be viewed as a suitable
``globalization'' of the closure of a single $N(\K)$--orbit in
$\Gr$. Moreover, $\ol{\Bun}_N$ has many advantages over $\Gr$ (see
\secref{expect}).

In this paper we define the appropriate ``Whittaker category'' of
perverse sheaves on $\ol{\Bun}_N$ (more precisely, on its
generalization $_{x,\infty}\BunNftb$) and prove that it does indeed
possess all the properties that one would expect from it by analogy
with the Whittaker functions. Namely, this category is semi-simple,
and each irreducible object of this category is a local system on a
single stratum extended by zero (these strata are the analogues of the
$N(\K)$--orbits in $\Gr$).

As an application of the main theorems of this paper, we compute the
cohomology of certain sheaves on the intersections of $G(\O)$-- and
$N(\K)$--orbits in $\Gr$. Using this result we prove our conjecture
from \cite{previous} in the case of a general reductive group $G$. In
the case of $G=GL_n$, B.C.~Ngo \cite{Ngo} has earlier given an elegant
proof of this conjecture using a different method.\footnote{In the
course of writing this paper we learned from B.C.~Ngo that he and
P.~Polo obtained an independent proof of the conjecture for a general
reductive group $G$.

{\em Note added in Nov 2000}: the paper by B.C.~Ngo and P.~Polo, {\em
R\'esolutions de Demazure affines et formule de Casselman-Shalika
g\'eom\'etriquehas}, has appeared as math.AG/0005022.} As explained in
\cite{previous}, the proof of this conjecture gives us a purely
geometric proof of the Casselman-Shalika formula \cite{CS,Shi} for the
Whittaker function.

More details on the background and motivations for the present work,
as well as the description of the structure of this paper, can be found
in \secref{intr}.

We note that the results of this paper are valid in the characteristic
$0$ case, when one works in the context of ${\mathcal D}$--modules.

\ssec*{Acknowledgements.} We are grateful to D. Kazhdan for valuable
discussions.

The research of E.F. has been supported by grants from the Packard
Foundation and the NSF. The research of D.G. has been supported by NSF
and K.V. has been receiving support from NSA and NSF.

\section{Background and overview}    \label{intr}

\ssec{Hecke algebra and Whittaker functions} \label{functback}

\sssec{} Let $G$ be a split connected reductive group over $\Fq$,
$G(\K)$ the corresponding group over the local field $\K=\Fq((t))$,
and $G(\O)$, its maximal compact group, where $\O = \Fq[[t]]$. We fix
a Haar measure on $G(\K)$ so that $G(\O)$ has measure $1$.

Consider the {\em Hecke algebra} ${\mathbf H}$ of $G(\K)$ with respect
to $G(\O)$, i.e., the algebra of $\Ql$--valued compactly supported
$G(\O)$--bi--invariant functions on $G(\K)$ with the convolution
product:

\begin{equation}    \label{conv1}
(h_1 \star h_2)(g) = \int_{G(\K)} h_1(x) h_2(g\cdot x^{-1}) \; dx.
\end{equation}

It is well--known that there is a bijection between the double quotient
of $G(\K)$ with respect to $G(\O)$ and the set $\Lambda^{++}$ of
dominant coweights of $G$: to $\la\in\Lambda^{++}$ we attach the coset
of $\la(t)\in T(\K)\subset G(\K)$, where $T$ is the Cartan subgroup of
$G$. Hence we obtain a basis $\{c_\la\}_{\la\in\Lambda^{++}}$ of
$\HH$ consisting of the corresponding characteristic functions.

It turns out that $\HH$ has a ``better'' basis. Let $\GL$ be the
Langlands dual group of $G$ and let $\on{Rep}(\GL)$ denote the
Grothendieck ring (over $\Ql$) of its finite-dimensional
representations.

The following statement is often referred to as the Satake isomorphism:

\sssec{\bf Theorem.}(\cite{Sa,La,G})   \label{Satake}
{\em There exists a canonical isomorphism of algebras $\HH \simeq
\on{Rep}(\GL)$.}

\sssec{}    \label{classala}
To an element $\la\in\Lambda^{++}$ we can attach the class of
the corresponding irreducible $\GL$--module $V^\la$ in $\on{Rep}(\GL)$.
Their images  under the Satake isomorphism, denoted by $A_\la$,  yield the
``better'' basis for $\HH$.

Moreover, the $A_\la$ have the following form:
\begin{equation}    \label{Ala}
A_\la=q^{-\langle\la,\rhoc\rangle}\left( c_\la + \sum_{\mu \in
\Lambda^{++};\, \mu<\la} p_{\la\mu}\cdot c_\mu \right), \quad \quad
p_{\la\mu} \in \Z[q].
\end{equation}

\sssec{}

Let $N\subset G$ be a maximal unipotent
subgroup. We define a character $\chi:N(\K)\to\Fq$ as follows:
$$\chi(u)=\sum_{i=1}^\ell \psi\left(\on{Res}(u_i)\right),$$ where
$u_i, i=1\ldots,\ell=\dim N/[N,N]$ are natural coordinates on
$N/[N,N]$ corresponding to the simple roots, $\on{Res}: \K \arr \Fq$ is the map
$\on{Res}\left( \sum_{n \in \Z} f_i t^i \right) = f_{-1}$ and $\psi:\Fq \arr
\Ql^{\times}$ is a fixed non-trivial character.

We define the space $\WW$ consisting of functions $f:G(\K)\to\Ql$,
such that

\begin{itemize}

\item
$f(g\cdot x)=f(g)$, if $x\in G(\O)$,

\item
$f(n\cdot g)=\chi(n)\cdot f(g)$, if $n\in N(\K)$,

\item
$f$ is compactly supported modulo $N(\K)$.

\end{itemize}

Clearly, the space $\WW$ is an $\HH$--module with respect to the
following action:
\begin{equation} \label{functheckeaction}
h\in\HH, f \in \WW \longrightarrow (f\star h)(g)= \int_{G(\K)}
f(g\cdot x^{-1})h(x)  \; dx.
\end{equation}
We call $\WW$ the {\em Whittaker module}.

There is a bijection between the double cosets $N(\K)\backslash
G(\K)/G(\O)$ and the coweight lattice $\Lambda$: $\la\to N(\K)\cdot \la(t)
\cdot G(\O)$.  However, it is easy to see that if $\la\in\Lambda$ is
non--dominant, then $f(\la(t))=0$ for any $f\in \WW$. Therefore $\WW$
has a basis $\{\phi_\la\}_{\la\in\lambda^{++}}$, where $\phi_\la$ is
the unique function in $\WW$ which is non--zero only on the double coset
$N(\K)\cdot\la(t) \cdot G(\O)$ and
$\phi(\la(t))=q^{-\langle\la,\rhoc\rangle}$
(here $\rhoc \in \check\Lambda$ is half the sum of positive roots of
$G$ and $\langle,\rangle$ is the canonical pairing of the weight and
coweight lattices).

\sssec{}
In \cite{previous} we considered a map $\FF:\HH\to\WW$ defined by the
formula $\FF(h)=\phi_0\star h$, i.e.,
$$(\FF(h))(g)=\int_{N(\K)} h(n^{-1}\cdot g) \chi(n)\; dn$$ (here the
measure on $N(\K)$ is chosen in such a way that the measure of $N(\O)$
is $1$). It is easy to see that $\FF$ is an isomorphism of
$\HH$--modules.

Remarkably, it turns out that $\FF$ is compatible with the above bases:

\sssec{\bf Theorem.}(\cite{previous})    \label{cs}
$\FF(A_\la)=\phi_\la$ {\em for all} $\la \in \Lambda^{++}$.

\sssec{Connection with the Whittaker functions} Each semi-simple
conjugacy class $\gamma$ of the group $\GL(\Ql)$ defines a
homomorphism $\ga: \on{Rep} \GL(\Ql) \arr \Ql$, which maps $[V]$
to $\on{Tr}(\gamma,V)$. We denote the corresponding homomorphism $\HH
\arr \Ql$ by the same symbol.

For each $\gamma$ as above one defines the {\em Whittaker function}
$W_\ga$ as the unique function on $G(\K)$ which satisfies:

\begin{itemize}
\item
$W_{\gamma}(g\cdot x) = W_{\gamma}(g), \forall x \in G(\O)$;

\item
$W_{\gamma}(n\cdot g) = \chi(n) \cdot W_{\gamma}(g), \forall\, n \in N(\K)$;

\item
$$h\star W_\ga = \gamma(h) \cdot W_{\gamma}, \quad \quad \forall
h \in \HH$$ (here $\star$ is defined by the same formula as in
\eqref{functheckeaction}).

\item
$W_\ga(1) = 1$.

\end{itemize}

It was proved in \cite{previous} that \thmref{cs} is equivalent to the
following well--known Casselman--Shalika formula (see \cite{CS,Shi}):
\begin{equation}    \label{cs1}
W_\gamma=\sum_{\la \in \Lambda^{++}} \on{Tr}(\gamma,(V^\la)^*) \cdot
\phi_\la.
\end{equation}

\bigskip

Thus, the functions $\phi_\la$ can be viewed as the ``building
blocks'' for the Whittaker function.

It is instructive to compare formula \eqref{cs1} to the formula for
the spherical function given in Sect.~7 of \cite{previous}; the
building blocks for the latter are the functions $A_\la$.

\sssec{} \thmref{cs} allows to compute $\phi_\mu\star A_\la$
for an arbitrary $\mu \in \Lambda^{++}$, as follows.

By \thmref{Satake}, the structure constants of the Hecke algebra are
equal to those of the Grothendieck ring $\on{Rep} \GL$ (the
``Clebsch--Gord\r{a}n coefficients''). We have:
$$A_\la \star A_\mu = \sum_{\nu \in \PL^+} C_{\la\mu}^\nu\, A_\nu,$$
where
$C_{\la\mu}^\nu = \dim \on{Hom}_{\GL}(V^\nu,V^\la \otimes V^\mu).$
Now \thmref{cs} implies that
\begin{equation} \label{act1}
\phi_\mu\star A_\la = \sum_{\nu \in \PL^+} C_{\la\mu}^\nu\, \phi_\nu.
\end{equation}

For $\mu\in\Lambda$, let $\chi_\mu:N(\K)\to\Ql^\times$ be a character
defined by the formula
$$\chi_\mu(n) = \chi(\mu(t) \cdot n \cdot \mu(t)^{-1}).$$

By evaluating both sides of \eqref{act1} on $g=(\la+\mu)(t)$
we obtain the following:

\begin{equation}    \label{eq2}
\int_{N(\K)}
A_{\la}(n^{-1}\cdot \nu(t))\, \chi_\mu(n)\, dn =
q^{-\langle \nu,\rhoch \rangle}\cdot C_{\la\mu}^{\mu+\nu}.
\end{equation}

\ssec{Geometrization}

\sssec{}

The goal of the present paper is to provide a geometric interpretation
of the results described in \secref{functback}, in particular,
of \thmref{cs} and formula \eqref{eq2}.

The starting point is the fact that the Hecke algebra $\HH$ admits a
natural geometric counterpart.

There exists an algebraic group (resp., an ind--group) over $\Fq$, whose
set of $\Fq$--points identifies with $G(\O)$ (resp, $G(\K)$); to
simplify notation, we will denote these objects by the same
symbols. The quotient $\Gr:=G(\K)/G(\O)$ is an ind--scheme and we can
consider the category $\on{P}_{G(\O_x)}(\Gr)$ of $G(\O)$--equivariant
perverse sheaves on $\Gr$, defined over the algebraic closure $\Fqb$
of $\Fq$ (since the closures of $G(\O)$--orbits in $\Gr$ are
finite-dimensional, there is no problem to define this category).

The $\Ql$--vector space $\HH$ is isomorphic to the Grothendieck group
of the category $\on{P}_{G(\O_x)}(\Gr)$. The irreducible objects of
$\on{P}_{G(\O_x)}(\Gr)$ are the intersection cohomology sheaves
$\A_\la$ attached to the corresponding $G(\O)$--orbits $\Gr^\la$. A
remarkable fact is that the class of $\A_\la$ in the Grothendieck
group corresponds (up to a sign) to the element $A_\la \in \HH$ defined in
\secref{classala}.

Moreover, \thmref{Satake} has a categorical version
(\thmref{sheaf-Satake}):
$\on{P}_{G(\O_x)}(\Gr)$ is a tensor category and as such it is
equivalent to the category ${{\mathcal R}ep}(\GL)$ of
finite--dimensional representations of $\GL$.

\sssec{}

Now we would like to find a geometric counterpart
$\on{P}_{N(\K)}^\chi(\Gr)$ of the Whittaker module $\WW$. A natural
candidate for it would be a category of $(N(\K),\chi)$--equivariant
perverse sheaves on $\Gr$, whose Grothendieck group would be
isomorphic to $\WW$.

Unfortunately, the situation here is much more complicated, since the
orbits of the group $N(\K)$ on $\Gr$ are infinite--dimensional and we
do not yet have a satisfactory theory of perverse sheaves with
infinite--dimensional supports.

In this paper we propose a substitute for the category of
$(N(\K),\chi)$--equivariant perverse sheaves on $\Gr$.  Our category
will consist of perverse sheaves, which, however, do not ``live'' on
$\Gr$, but rather on an algebraic stack $\ol{\Bun}_N$, which is
defined using a global curve $X$ (see \secref{right category}).

\sssec{}    \label{expect}

The stack $\ol{\Bun}_N$ is the Drinfeld compactification of the moduli
space of $N$--bundles on $X$. This stack has been used in dealing with 
several other problems of geometric representation theory \cite{FM,FM2,BG1}.

Let us sketch the definition in the simplest example of $G=GL(2)$. In this 
case, the stack $\Bun_N$ classifies short exact sequences
$$0\to \OO_X\to \E\to \OO_X\to 0,$$ where $\OO_X$ is the structure
sheaf on $X$ and $\E$ is a rank two vector bundle. In other words, it
classifies rank two vector bundles $\E$ with $\on{det} \E \simeq
\OO_X$, together with an embedding of the trivial line bundle into
$\E$ (i.e., a {\em maximal} embedding of sheaves $\OO_X
\hookrightarrow \E$).

The ``compactification'' $\ol{\Bun}_N$ of $\Bun_N$ classifies rank two
vector bundles $\E$ (with $\on{det} \E \simeq \OO_X$), together with an
{\em arbitrary} embedding of sheaves $\OO_X \hookrightarrow \E$. 

The definition in the general case follows the same lines and we refer
the reader to \secref{compact} for details.

\medskip

Informally, from the point of view of $\on{P}_{N(\K)}^\chi(\Gr)$, the
idea of introducing $\ol{\Bun}_N$ is as follows.  Since
$N(\K)$--orbits on $\Gr$ are infinite--dimensional, it is natural to
try to ``quotient them out'' by a subgroup of $N(\K)$, which acts
almost freely on the $N(\K)$--orbits in $\Gr$. There is a natural
candidate for such a subgroup: fix a projective curve $X$ and a point
$x\in X$ (so that $\K\simeq \K_x$) and consider the group
$N_{\out}\subset N(\K)$ of maps $(X\setminus x) \ \to \ N$.

The quotient $S/N_{\out}$ of a single $N(\K)$--orbit $S$ on $\Gr$
makes perfect sense: this is the stack $\Bun_N$. In
\secref{introrbits} we give a precise scheme-theoretic definition of
the closure $\ol{S}$ of $S$. This is an ind-scheme, which appears to
be highly non-reduced (see \remref{nilpotents}). Unlike the case of
$S$, it is not quite clear what the quotient of $\ol{S}$ by $N_{\out}$
should be, since in general the quotient of an ind-scheme by an
ind-group is not an algebraic stack.

In fact, there is a natural morphism $\ol{S} \to \ol{\Bun}_N$. The
image of the corresponding map at the level of $\Fq$--points is the
set-theoretic quotient $\ol{S}/N_{\out}$. But in general it is a
constructible subset of $\ol{\Bun}_N(\Fq)$. For example, in the case
of $G=SL_2$, locally this subset looks like $\AA^2\setminus
(\AA^1\setminus \on{pt})$ inside $\AA^2$. Thus, there is no obvious
way to interpret the naive set-theoretic quotient $\ol{S}/N_{\out}$ as
the set of points of an algebraic variety or an algebraic stack. On
the other hand, $\ol{\Bun}_N$ is a well-defined algebraic stack, and
so it may be viewed as an appropriate replacement for
$\ol{S}/N_{\out}$.

An important feature of $\ol{\Bun}_N$ is that it has a stratification
analogous to the stratification of $\ol{S}$ by $N(\K)$--orbits (see
\corref{stratified}). But the dimensions of the strata on $\ol\Bun_N$
are all even and equal to twice the naive (relative) dimensions of the
corresponding $N(\K)$--orbits on $\Gr$.

\sssec{}  \label{expprop}

The stack $\ol{\Bun}_N$ (more precisely, its generalization
${}_{x,\infty}\BunNftb$ introduced in \secref{variants}) allows us to
introduce the sought after category $\on{P}_{N(\K)}^\chi(\Gr)$. It
will possess the following basic properties, motivated by the
analogous properties of the Whittaker module $\WW$:

\begin{itemize}

\item[(1)]
The tensor category $\on{P}_{G(\O_x)}(\Gr)$ acts on
$\on{P}_{N(\K)}^\chi(\Gr)$ by ``Hecke convolution functors''.

\item[(2)]
The irreducible objects $\{\Psib^\la\}$ in $\on{P}_{N(\K)}^\chi(\Gr)$ are
labeled by $\la\in \Lambda^{++}$ and the functor
$\F:\on{P}_{G(\O_x)}(\Gr)\to \on{P}_{N(\K)}^\chi(\Gr)$ given by
$\F(\S)\to \Psib^0\star \S$ sends $\A_\la$ to $\Psib^\la$.

\item[(3)]
Each $\Psib^\la$ (in spite of being an irreducible perverse sheaf) is
an extension by $0$ of a local system from the corresponding stratum
of the stack $\BunNftb$.

\end{itemize}

\sssec{} As we will see in \secref{proofofsemisimplicity}, the above
properties (1) and (2) of $\on{P}_{N(\K)}^\chi(\Gr)$ formally imply
that this category is semi--simple and that the functor $\F$ is an
equivalence of categories. This will in turn
imply (almost tautologically) property (3).

Property (3) of the sheaves $\Psib^\la$ (which are the sheaf theoretic
counterparts of the functions $\phi_\la$) is perhaps the most
intriguing aspect of the category $\on{P}_{N(\K)}^\chi(\Gr)$ and is
the key point of this paper. We repeat that it
means that $\Psib^\la$ has zero stalks
outside of the locus where it is a local system. An irreducible
perverse sheaf satisfying this property is called ``clean'' and the
appearance of such perverse sheaves in representation theory always
has remarkable consequences.

To contrast this to the case of $G(\O)$--equivariant sheaves, note
that while the category of $\on{P}_{G(\O)}(\Gr)$ is semi-simple, the
irreducible objects ${\mc A}_\la$ of $\on{P}_{G(\O)}(\Gr)$ have
non-zero stalks on all $G(\O)$--orbits that lie in the closure of the
orbit $\Gr^\la$.

\sssec{Remark.} As we have seen above, $N(\K)$--orbits on $\Gr$ are
enumerated by $\la \in \Lambda$ but only those corresponding to $\la
\in \Lambda^{++}$ can carry elements of $\WW$.

However, the stack $_{x,\infty}\BunNftb$ will have many more strata,
among which we single out what we call the relevant ones (the latter
will be in bijection with $\Lambda^{++}$). The perverse sheaves
belonging to $\on{P}_{N(\K)}^\chi(\Gr)$ will have non-zero stalks only
on the relevant strata.

\sssec{}

Another source of motivation for us was the following.
Although our definition of
$\on{P}_{N(\K)}^\chi(\Gr)$ is not obvious, it is easy to formulate a
sheaf counterpart of formula \eqref{eq2}.

Denote by $\ol{\Gr}^\la$ the closure of the orbit $G(\O)\cdot \la(t)$
in $\Gr$ and by $S^\nu$ the $N(\K)$--orbit of $\nu(t)\in\Gr$. The
variety $\ol{\Gr}^\la\cap S^\nu$ is finite--dimensional and if
$\mu\in\Lambda$ is such that $\mu+\nu$ is dominant, the character
$\chi_\mu:N(\K)\to\GG_a$ gives rise to a map
$\chi_\mu^\nu:\ol{\Gr}^\la\cap S^\nu\to\GG_a$
(cf. \secref{admcharacters} for the definition of $\chi_\mu^\nu$).
Let $\J_\psi$ be the Artin--Schreier sheaf on $\GG_a$.

The following is a generalization of our Conjecture 7.2. from
\cite{previous} (in \cite{previous} it was stated for $\mu=0$, and in
that form it was subsequently proved by B.~C.~Ngo \cite{Ngo} for
$G=GL_n$). 

\begin{thmm}    \label{gencs1}
For $\la\in\Lambda^{++}$ and $\mu,\nu\in\Lambda$ with
$\mu+\nu\in\Lambda^{++}$ the cohomology
\begin{equation} \label{gencs}
H_c^k(\ol{\Gr}^\la\cap S^{\nu}, \A_\la|_{\ol{\Gr}^\la\cap S^\nu}\otimes
\chi^\nu_{\mu}|_{\ol{\Gr}^\la\cap S^\nu}{}^*(\J_\psi))
\end{equation}
vanishes unless $k=\langle 2\nu,\rhoc\rangle$ and
$\mu\in\Lambda^{++}$. In the latter case, this cohomology identifies
canonically with $\Hom_{\GL}(V^{\la}\otimes V^{\mu},V^{\mu+\nu})$.
\end{thmm}

At the end of this paper we will prove this theorem using our category
$\on{P}_{N(\K)}^\chi(\Gr)$ and mainly the ``cleanness'' property of
the objects $\Psib^\la$ mentioned above.\footnote{This approach to the
Casselman-Shalika formula was suggested to one of us (D.G.) by his thesis
advisor, J.~Bernstein, several years ago.}

It will turn out (as in the function--theoretic setting) that the
stalks of the convolution $\Psib^{\mu+\nu}\star\A_{-w_0(\la)}$ are
isomorphic to the cohomology groups appearing in \eqref{gencs}.

By passing to the traces of the Frobenius, the above theorem entails
formula \eqref{eq2}, and hence (in the  case $\mu=0$) the
Casselman--Shalika formula \eqref{cs1}, as explained in Sect.~6 of
\cite{previous}. Thus, we obtain a geometric proof of the
Casselman-Shalika formula.

\ssec{Contents}

The main part of the paper is devoted to the study of the Drinfeld
compactification $\BunNftb$ and a certain category of perverse sheaves
on it.


We start out in \secref{compact} with the definition of the stack
$\BunNftb$ and its ind--version $_{\xl,\infty}\BunNftb$.


In \secref{pproperties} we focus on the properties of $\BunNftb$ that
are related to the $N(\K)$--action. An important fact proved in
\secref{pproperties} is that the embedding
$\BunNft\hookrightarrow\BunNftb$ is affine, which is an analogue of
the fact that the embedding $S^\nu\hookrightarrow \ol{S}^\nu$ is
affine. In addition, we prove that $\BunNft$ is topologically
``contractible'' (by analogy with $S^\nu$, which is isomorphic to an
infinite--dimensional affine space).


In \secref{category} we introduce our main objects of study--the
perverse sheaves $\Psib_\varpi$ on $\ol{\Bun}_N$ or, rather, on its
twisted version $\BunNftb$. We formulate \thmmref{oneprime}, which expresses
the cleanness property of the perverse sheaves $\Psib_\varpi$, and its generalization,
\thmmref{one}, which is the main result of this paper. 


In \secref{hecke} we introduce the main tool needed for the proof of
\thmmref{one} -- the action of the category $\on{P}_{G(\O_x)}(\Gr)$ on
the derived category of sheaves on $_{x,\infty}\BunNftb$ by Hecke
functors.


In \secref{proofs} we derive \thmmref{one} from \thmmref{two}, which
describes the action of the Hecke functors on our basic perverse
sheaves.


\secref{heckecontd} is devoted to the proof of \thmmref{two}.


Finally, in \secref{conj} we prove a generalization of the conjecture
from \cite{previous} (see \thmmref{gencs1} above).

\ssec{Notation and conventions} \label{notation}

{}From now on, we replace the field $\Fq$ by its algebraic closure
$\Fqb$.


Throughout the paper, $X$ will be a fixed smooth projective connected
curve over $\Fqb$.


Furthermore, $G$ will stand for a connected reductive group over
$\Fqb$. We fix a Borel subgroup $B\subset G$ and let $N$ be its
unipotent radical and $T = B/N$. We denote by $\Lambda$ the {\it
coweight} lattice and by $\check\Lambda$ the {\it weight} lattice of
$T$; $\langle,\rangle$ denotes the natural pairing between the two.

The set of vertices of the Dynkin diagram of $G$ is denoted by $\I$;
$\Lambda^{++}$ denotes the semigroup of dominant coweights and
$\Lambda^+$ the span over $\Z_{+}$ of positive coroots;
$\check\Lambda^{++}$ and $\check\Lambda^{+}$ denote similar objects
for $\check\Lambda$; $\rhoc\in\check\Lambda$ is the half sum of
positive roots of $G$ and $w_0$ is the longest element of the Weyl
group.

For $\check\la\in\check\Lambda^{++}$, we write $\V^{\check\la}$ for
the corresponding Weyl module, i.e.,
$$\V^{\check\la}=\Gamma(G/B,\L_{-w_0(\check\lambda)})^*,$$ where
$\L_{-w_0(\check\lambda)}$ is the line bundle on $G/B$ corresponding to
the character $-w_0(\check\lambda)$ of $T$. We obtain a collection of
canonical embeddings
$\V^{\lambdach+\much} \arr
\V^{\lambdach} \otimes \V^{\much}$ for all $\lambdach,\much \in
\check\Lambda^{++}$.

For $\la \in \Lambda^{++}$, we denote by $V^\la$ the irreducible
representation of $\GL(\Ql)$ with highest weight $\la$.

\medskip

This paper will extensively use the language of algebraic stacks (over
$\Fqb$), see \cite{La:s}. When we say that a stack $\Y$ classifies
{\it something}, it should always be clear what an $S$--family of {\it
something} is for any $\Fqb$--scheme $S$, i.e., what is the groupoid
$\Hom(S,\Y)$, and what are the functors $\Hom(S_2,\Y) \arr
\Hom(S_1,\Y)$ for each morphism $S_1 \arr S_2$.

For example, if $H$ is an algebraic group, we define $\Bun_H$ as a
stack that classifies $H$--bundles on $X$. This means that
$\on{Hom}(S,\Bun_H)$ is the groupoid whose objects are $H$--bundles on
$X \times S$ and morphisms are isomorphisms of these
bundles. Pull--back functor for $S_1\to S_2$ is defined in a natural
way.

\medskip

Let $x\in X$ be a point.  We will denote by $\O_x$ (resp., $\K_x$) the
completion of $\OO_X$ at $x$ (resp., the field of fractions of
$\O_x$).  We will use the notation $\D_x$ (resp., $\D_x^\times$) for
the formal (resp., formal punctured) disc around $x$ in $X$. They
will not appear as schemes and we will use them only in the following
circumstances:

\noindent Let $S$ be a scheme and assume for simplicity that $S$ is
affine, $S=\on{Spec}(\OO_S)$.  An $S$--family of vector bundles on
$\D_x$ (resp., $\D^\times_x$) is by definition a projective finitely
generated module over the completed tensor product
$\OO_S\widehat\otimes \O_x$ (resp., $\OO_S\widehat\otimes
\K_x$). Similarly, one defines an $S$--family of $G$--bundles on
$\D_x$ or $\D_x^\times$.

\noindent Let $S$ be as above and let $Y$ be another affine scheme
$Y=\on{Spec}(\OO_Y)$. An $S$--family of maps $\D_x\to Y$ (resp.,
$\D^\times_x\to Y$) is by definition a ring homomorphism $\OO_Y\to
\OO_S\widehat\otimes \O_x$ (resp., $\OO_Y\to \OO_S\widehat\otimes
\K_x$).

\medskip

This paper deals with perverse sheaves on algebraic stacks. Although
these objects make a perfect sense, the corresponding derived category
is problematic: there is no doubt that the ``correct'' derived
category exists, but a satisfactory definition is still unavailable in
the published literature.

For the needs of our paper, we could stay within the abelian category of
perverse sheaves, and formulate all of our statements and proofs in
terms of complexes of perverse sheaves and their
cohomologies. However, doing so would be very inconvenient as it would
force us to use cumbersome notation and increase the length of the
proofs. For this reason, in this paper we will adopt the point of view
that the appropriate definition of the derived category $\Sh(\Y)$ on a
stack $\Y$ exists. The interested reader can easily reformulate all of the
statements below so as to avoid the use of the derived category.

Thus, when we discuss objects of the derived category, the
cohomological gradation should always be understood in the perverse
sense. In addition, for a morphism $f:\Y_1\to \Y_2$, the functors
$f_{!}$, $f_{*}$, $f^*$ and $f^{!}$ should be understood ``in
the derived sense''.

\medskip

Let $\Y$ be a stack, $H$ an algebraic group and $\F_H$ an $H$--torsor
over $\Y$. Let, in addition, $Z$ be an $H$--scheme. We denote by
$\F_H\underset{H}\times Z$ the corresponding fibration over $\Y$, with
the typical fiber $Z$. (Unfortunately, the notation
$A\underset{C}\times B$ is also used for the fiber product of $A$ and
$B$ over $C$, but in all instances it will be clear from the context
what the notation means).

Now, if $\T$ is a perverse sheaf (or a complex) on $\Y$ and $\S$ is an
$H$--equivariant perverse sheaf (or complex) on $Z$, we can form their
twisted external product and obtain a perverse sheaf (resp., a
complex) on $\F_H\underset{H}\times Z$, denoted $\T\tboxtimes\S$,
which is $\T$ ``along the base'' and $\S$ ``along the fiber''.

\section{Drinfeld's compactification}  \label{compact}

{}From now on, with the exception of \secref{conj}, we will work under
the assumption that the derived group $[G,G]$ is {\em
simply-connected}.

\ssec{The stack $\Bun_N$}

\sssec{}

Consider the moduli stack $\Bun_G$ of principal $G$--bundles over the
curve $X$. We repeat that
by definition, for an $\Fqb$--scheme $S$, $\on{Hom}(S,\Bun_G)$ is the
groupoid, whose objects are $G$--bundles on $X \times S$ and morphisms
are isomorphisms of these bundles. One defines similarly the stacks
$\Bun_B$, $\Bun_T$ and $\Bun_N$.

We have natural morphisms of stacks
$$\po:\Bun_B\to\Bun_G, \quad \quad {\mathfrak q}:\Bun_B\to\Bun_T,$$
that send a given $B$--bundle $\F_B$ to the bundles
$$\po(\F_B)=\F_G:=\F_B\underset{B}{\times} G, \quad \quad
\qo(\F_B)=\F_T:=\F_B\underset{B}\times T=\F_B/N,$$ respectively.

The projection $\po$ is representable, but in general is not smooth
and has non-compact fibers, while the projection $\qo$ is smooth, but
in general non-representable.

Denote by $\F^0_T$ the trivial $T$--bundle on $X$. It follows from the
definitions that the closed substack $\qo^{-1}(\F^0_T)\subset\Bun_B$
can be identified with the stack $\Bun_N$. In
what follows, for a fixed $T$-bundle $\F_T$ we will denote by
$\BunNft$ the closed substack $\qo^{-1}(\F_T) \subset \Bun_B$.

Here is a Pl\"ucker type description of the stack $\BunNft$:

\sssec{} \label{def-bunN} Given a scheme $S$, $\on{Hom}(S,\bn)$ is the
groupoid whose objects are pairs $(\F_G,\map)$, where $\F_G$ is a
$G$--bundle on $X \times S$ and $\map = \{ \map^{\lambdach} \}$ is a
collection of {\it maximal} embeddings
\begin{equation}
\map^{\lambdach}:\L_{\F_T}^{\lambdach}\hookrightarrow
\V^{\lambdach}_{\F_G},\quad \quad \forall\lambdach\in\Lambdach^{++}.
\end{equation}
Here $\L_{\F_T}^{\lambdach}$ is the line bundle on $X \times S$
induced from $\F_T$ by means of the character $\lambdach:T\to{\mathbb
G}_m$, $\V^{\lambdach}_{\F_G}$ is the vector bundle on $X \times S$
associated with the representation $\V^{\lambdach}$ of $G$ and the
$G$--bundle $\F_G$.

Recall that an embedding of locally free sheaves is called maximal if
it is a bundle map, i.e., an injective map of coherent sheaves,
such that the quotient is torsion-free.

The above collection $\map$ of embeddings must satisfy the so-called
Pl\"ucker relations. Namely, given a pair of dominant integral weights
of $G$, $\lambdach$ and $\much$, the map
$$\L_{\F_T}^{\lambdach}\otimes
\L_{\F_T}^{\much}\overset{\map^{\lambdach}\otimes\map^{\much}}
\longrightarrow \V^{\lambdach}_{\F_G}\otimes \V^{\much}_{\F_G}\simeq
(\V^{\lambdach}\otimes \V^{\much})_{\F_G}$$ must coincide with the
composition
$$\L_{\F_T}^{\lambdach}\otimes \L_{\F_T}^{\much}\simeq
\L_{\F_T}^{\lambdach+\much}
\overset{\map^{\lambdach+\much}}\longrightarrow
\V^{\lambdach+\much}_{\F_G}\hookrightarrow (\V^{\lambdach}\otimes
\V^{\much})_{\F_G},$$ where the last arrow comes from the canonical
embedding of representations $\V^{\lambdach+\much}\to
\V^{\lambdach}\otimes \V^{\much}$ (see \secref{notation}).

The morphisms between two objects $(\F_G,\map)$ and $(\F'_G,\map')$ in
$\on{Hom}(S,\bn)$ are isomorphisms between $\F_G$ and $\F'_G$, which
render the corresponding diagrams involving the maps
$\map^{\lambdach}, \map'{}^{\lambdach}$ commutative.

The definition of the functor $\on{Hom}(S_2,\bn) \arr
\on{Hom}(S_1,\bn)$ associated to a morphism $S_1 \arr S_2$ is
straightforward.

\sssec{Example of $G=GL_2$} \label{example}

Given a line bundle $\L$ on $X$, let us denote by $\E_{\L}$ the
algebraic stack classifying the short exact sequences:
\begin{equation}    \label{exten}
0\to \L\to \E \to \OO_X \to 0.
\end{equation}

More precisely, the objects of the groupoid $\on{Hom}(S,\E_L)$ are
coherent sheaves $\E$ on $X \times S$ together with a short exact
sequence $$0 \arr \L \boxtimes \OO_S \arr \E \arr \OO_X \boxtimes
\OO_S \arr 0,$$ and morphisms are morphisms between such exact sequences
which are identities at the ends.

There exists a canonical map from $\E_{\L}$ to the affine space
$H^1(X,\L)$ which associates to every short exact sequence as above
its extension class. Moreover, it is easy to see that the stack $\E_L$
is isomorphic to the quotient of $H^1(X,\L)$ by the trivial action of
the additive group $H^0(X,\L)$.

Now, for $G=GL_2$, a $T$--bundle $\F_T$ is the same as a pair of line
bundles $\L_1$ and $\L_2$. Set $\L=\L_1\otimes\L_2^{-1}$. According to
the definition, the stack $\BunNft$ is naturally isomorphic to
$\E_{\L}$.

\ssec{Compactification}  \label{defcompact}

\sssec{}

As we remarked above, the fibers of the projection
$\po:\Bun_B\to\Bun_G$ are non-compact. The reason is that the variety
of line subbundles of a fixed degree in a given vector bundle is
non-compact. However, the latter has a natural compactification: the
variety of invertible subsheaves of a fixed degree in a vector bundle,
considered as a locally free sheaf. Following this example,
V.~Drinfeld proposed the following (partial) compactification of the stack
$\BunNft$ along the fibers of the projection $\po$.

\sssec{\bf Definition.}
The stack $\BunNftb$ classifies pairs
$(\F_G,\map)$, where $\F_G$ is a $G$--bundle on $X$, and $\map = \{
\map^{\lambdach}\}$ is the collection of maps
$$\map^{\lambdach}:\L_{\F_T}^{\lambdach}\hookrightarrow
\V^{\lambdach}_{\F_G},$$ which are now embeddings of coherent sheaves
(i.e., we do not require any more that the quotient is
torsion-free). These embeddings must satisfy the Pl\"ucker relations
in the same sense as in Definition \ref{def-bunN}.

More precisely, the objects of the groupoid $\on{Hom}(S,\BunNftb)$ are
pairs $(\F_G,\map)$, where $\F_G$ is a $G$--bundle on $X \times S$ and
$\map = \{ \map^{\lambdach} \}$ is a collection of embeddings
\begin{equation}
\map^{\lambdach}: \L_{\F_T}^{\lambdach}\hookrightarrow
\V^{\lambdach}_{\F_G},\quad \quad \forall\lambdach\in\Lambdach^{++},
\end{equation}
such that the quotient $\V^{\lambdach}_{\F_G}/\L_{\F_T}^{\lambdach}$ is
$S$--flat and the Pl\"ucker relations hold.
Morphisms in $\on{Hom}(S,\BunNftb)$ are defined in the
same way as in $\on{Hom}(S,\BunNft)$.

\sssec{} It is clear that $\BunNft$ is an open substack of $\BunNftb$,
and we denote by $j$ the corresponding open embedding.
One can show that $\BunNft$ is dense in $\BunNftb$. We will not
need this fact and the reader is referred to \cite{BG} for the proof.

We will denote by $\p:\BunNftb\to\Bun_G$ the natural
projection:
$$\p(\F_G,\map)=\F_G.$$

It is clear that the morphism $\p$ extends $\po$, i.e., that
$\po=\p\circ j$. Moreover, we have the following statement, which
follows from the fact that the variety of invertible subsheaves of a
fixed degree in a given vector bundle is a complete variety.

Consider the natural $T$--action on $\BunNftb$: a point $\tau\in T$
acts on a point $(\F_G,\map)\in\BunNftb$ by multiplying each
$\map^{\lambdach}$ by $\lambdach(\tau)$. The projection $\p$ clearly
factors as $\BunNftb\to \BunNftb/T\to \Bun_G$.

\sssec{\bf Lemma.}
{\em The morphism $\p:\BunNftb/T\to\Bun_G$ is representable and proper.}

\bigskip

This is the reason why $\BunNftb$ is called a ``compactification'' of
$\bn$.

\sssec{Example.} Let us again consider the case of $G=GL_2$ using
the notation of \secref{example}.

In this case $\BunNftb$ is the stack that classifies the triples
$(\E,\map_1,\map_2)$, where $\E$ is a rank $2$ bundle on $X$, $\map_1$
(resp., $\map_2$) is a non-zero map $\L_1\to\E$ (resp., $\E\to\L_2$)
such that the composition $\map_2\circ\map_1$ vanishes and the induced
map
$$
\on{det}(\E)\to \L_1\otimes\L_2
$$
is an isomorphism.

The fiber of $\p$ over a rank $2$ bundle $\E$ is the vector space
$\on{Hom}(\L_1,\E)$ with the origin removed.

Note that $\BunNftb$ is stratified by locally closed substacks
$_d \BunNftb, d\geq 0$, which classify those triples
$(\E,\map_1,\map_2)$, for which the divisor of zeros of $\map_1$ is
of order $d$. In particular, $_0 \BunNftb = \BunNft$.

Let us describe the $\Fqb$--points of the substacks $_d \BunNftb$. A
point of $\BunNft$ is an isomorphism class of extensions
$$
0\to \L_1 \to \E \to \L_2 \to 0.
$$
Hence the set of $\Fqb$--points of $\BunNft$ is $H^1(X,\L)$, where
$\L=\L_1\otimes\L_2^{-1}$. A point of $_1 \BunNftb$ is an isomorphism
class of extensions
$$
0 \to \L_1(x) \to \E \to \L_2(-x) \to 0,
$$
with some fixed $x \in X$. Therefore the set of $\Fqb$--points of $_1
\BunNftb$ is the set of points of a vector bundle over $X$ with fiber
$H^1(X,\L(2x))$ over $x \in X$. Likewise, we find that the set of
$\Fqb$--points of $_d \BunNftb$ is in one-to-one correspondence with
the set of points of a vector bundle over $S^d X$ with the fiber
$H^1(X,\L(2x_1 + \ldots + 2 x_d))$ over the point $x_1 + \ldots + x_d
\in S^i X$.

\sssec{} The following statement shows that we have a similar
stratification of $\BunNftb$ for general $G$. It is here that the
assumption that $[G,G]$ is simply-connected, which was made at the
beginning of this section, becomes essential.

Suppose that $\F'_T$ is a $T$--bundle, such that
$$
\L_{\F'_T}^{\lambdach} = \L_{\F_T}^{\lambdach}\left( \sum_k \langle
\nu_k,\lambdach \rangle \cdot x_k \right), \quad \quad \forall
\lambdach \in \Lambdach.
$$
Then we will write:
$$
\F'_T = \F_T \left( \sum_k \nu_k \cdot x_k \right).
$$

\sssec{\bf Proposition.}    \label{descrBunNftb}
{\em Let $(\F_G,\map)$ be an
$\Fqb$--point of $\BunNftb$. Then there exists a unique divisor $D =
\sum_k \nu_k\cdot x_k$ with $\nu_k\in \Lambda^{+}$ such that for
$\F'_T = \F_T(D)$ the meromorphic maps
$$\L^{\lambdach}_{\F'_T}\simeq \L^{\lambdach}_{\F_T}\left( \sum_k
\langle\nu_k,\lambdach\rangle\cdot x_k \right) \to
\V^{\lambdach}_{\F_G}$$
are regular everywhere and maximal.}

\sssec{Proof.}
For each $\i\in \I$, pick a fundamental weight
$\check\omega_\i\in\check\Lambda^{++}$.

Let $D_\i$ be the divisor of zeros of the map
$$\map^{\omegach_\i}:\L_{\F_T}^{\omegach_\i}\to \V^{\omegach_\i}_{\F_G}.$$
Set $D=\sum_{\i\in\I} D_{\i} \cdot \al_{\i}$. It is easy to check that $D$
satisfies all the requirements.

\sssec{\bf Corollary.}    \label{stratified}

{\em $\BunNftb$ is stratified by locally closed substacks
$_{\gamma}\Bun_N^{\F_T}$, with $\gamma = - \sum_{\i\in\I} d_\i \al_\i
\in -\Lambda^+$. The substack $_{\gamma}\Bun_N^{\F_T}$ is a bundle
over $\prod_{\i\in\I} S^{d_\i} X$ whose fiber
$_{\ol{x},\gamma}\Bun_N^{\F_T}$ at
$$
\ol{x} = (x_{\i,1} + \ldots + x_{\i,d_\i})_{\i\in\I}
$$
is isomorphic to $\Bun_N^{\F'_T}$, where} $$\F'_T = \F_T\left(
\sum_{\i\in\I} \sum_{m_\i=1}^{d_\i} \al_\i \cdot x_{\i,m_\i}
\right).$$

\bigskip

The above substacks $_{\ol{x},\gamma}\Bun_N^{\F_T}$ give us a
stratification of $\BunNftb$. We note that the codimension of the
stratum $_{\ol{x},\gamma}\Bun_N^{\F_T}$ in $\BunNftb$ equals
$2\sum_{\i\in\I} d_\i$; thus, it is always even.

\sssec{Remark.}

Consider the stack $\BunNftb/T$ for $G=GL(2)$. In \cite{B} and
\cite{Th} it was shown that an appropriate semi-stability condition
defines in $\BunNftb/T$ an open substack, which is in fact an
algebraic variety, and which can be constructed by a series of
blow-ups and blow-downs in the projective space ${\mathbb P}
H^1(X,\L)$.

Moreover, this description (referred to as ``geometric
approximation'') defines on this open substack of $\BunNftb/T$ a
stratification that coincides with that of \corref{stratified}.

For general $G$, such a description of a semi-stable part of
$\BunNftb/T$ appears to be unknown, and it would be interesting to
find one.

\ssec{Generalizations.} \label{variants}

\sssec{}

In addition to the stack $\BunNftb$, it will be convenient for us to
consider its generalization described below.

Let $\xl=\{x_1,...,x_n\}$ be a collection of $n$ distinct points on
$X$ and let $\nul$ be an $n$-tuple of elements of $\Lambda$.  Let
$_{\xl,\nul}\BunNftb$ be the stack classifying the data of
$(\F_G,\map)$, where $\map = \{ \map^{\lambdach} \}$ is now a
collection of arbitrary non-zero maps of coherent sheaves
$$\map^{\lambdach}:\L_{\F_T}^{\lambdach}\to \V^{\lambdach}_{\F_G}\left(
\sum_k \langle\nu_k,\lambdach\rangle\cdot x_k \right),$$
subject to the Pl\"ucker relations.

In particular, when $n=1$, $\xl=\{x\}$, $\nul=\{\nu\}$, we will write
$_{x,\nu}\BunNft$ instead of $_{\xl,\nul}\BunNftb$.

By definition, the $\BunNftb$ is nothing but $_{\xl,\zerol}\BunNftb$
and for a general $\nul$ we have a natural isomorphism
\begin{equation}    \label{natiso}
_{\xl,\nul}\BunNftb\simeq \overline{\Bun}_N^{\F'_T},
\end{equation}
where
$$
\F'_T = \F_T \left( -  \sum_k \nu_k \cdot x_k \right).
$$

We will use the notation $\nul'\geq\nul$ if for every
$k=1,...,n$, $\nu'_k-\nu_k\in\Lambda^+$. It is clear that if
$\nul'\geq\nul$, then there is a natural closed embedding
$$_{\xl,\nul}\BunNftb\hookrightarrow {}_{\xl,\nul'}\BunNftb.$$

We define the ind-stack $_{\xl,\infty}\BunNftb$ to be the inductive
limit of ${_{\xl,\nul}\BunNftb}$ as $\nul \in (\Lambda^+)^n$.

Note that for $\nul,\nul'\in\Lambda^n$, the closed substacks
$_{\xl,\nul}\BunNftb$ and $_{\xl,\nul'}\BunNftb$ belong to the same
connected component of $_{\xl,\infty}\BunNftb$ only if for every $k$
the projection of $\nu'_k-\nu_k$ to
$\Lambda/\on{Span}\{\Lambda^+\}\simeq \pi_1(G)$ equals $0$.

In what follows, we will continue to denote by $\p$ the natural map
$_{x,\infty}\BunNftb\to\Bun_G$.

\sssec{}

We will also need the following locally closed substacks of
$_{\xl,\infty}\BunNftb$.

Let us denote by
$$j_{\nul}: {}_{\xl,\nul} \Bun{}^{\F_T}_N
\hookrightarrow {}_{\xl,\infty}\BunNftb$$ the locally closed substack
corresponding to those pairs $(\F_G,\map)$, for which each map
$$\map^{\lambdach}: \L_{\F_T}^{\lambdach} \to
\V^{\lambdach}_{\F_G}\left( \sum_k \langle\nu_k,\lambdach\rangle \cdot
x_k \right)$$ is regular and maximal for all
$\lambdach\in\check\Lambda^{++}$ on the whole of $X$.

Note that under the isomorphism \eqref{natiso}, $_{\xl,\nul}
\Bun{}^{\F_T}_N$ goes to $\Bun^{\F'_T}_N$.

Let us also introduce the locally-closed substack
$$\wt{j}_{\nul}:{}_{\xl,\nul}\wt{\Bun}^{\F_T}_N \hookrightarrow
{}_{\xl,\infty}\BunNftb$$ as the locus corresponding to pairs
$(\F_G,\map^{\lambdach})$ for which each map
$$\map^{\lambdach}:\L_{\F_T}^{\lambdach}\to \V^{\lambdach}_{\F_G}\left(
\sum_k \langle \nu_k,\lambdach \rangle\cdot x_k \right)$$ is regular
and is moreover maximal {\it in a neighborhood of
$\underset{k=1}{\overset{n}\bigcup} x_k$} for all
$\lambdach\in\check\Lambda^{++}$. Note that this is equivalent to
saying that the quotient $\V^{\lambdach}_{\F_G}\left( \sum_k \langle
\nu_k,\lambdach \rangle\cdot x_k \right)/\L_{\F_T}^{\lambdach}$ has no
torsion supported at any of the points $x_1,...,x_n$.

To summarize, for each $\nul$ we have a sequence of embeddings:
$$_{\xl,\nul}\BunNft \hookrightarrow
{}_{\xl,\nul}\wt{\Bun}{}^{\F_T}_N \hookrightarrow
{}_{\xl,\nul}\BunNftb\hookrightarrow {}_{\xl,\infty}\BunNftb,$$ where
the first two arrows are open embeddings and the last arrow is a
closed one.

\section{Properties of $\BunNftb$.}  \label{pproperties}

In this section we will prove several important technical facts
concerning the structure of $\BunNftb$. These facts
will not be used until \secref{proofs}.

\ssec{Presentation as a double quotient.}    \label{dblquot}

\sssec{The structure of $N(\K_y)$.} Let $y$ be an $\Fqb$--point of $X$.
Let $G(\O_y)$ (resp., $G(\K_y)$) be be the group scheme (resp.,
group ind--scheme) that classifies maps $\D_y\to G$ (resp.,
$\D_y^\times\to G$) (cf. \secref{notation}). In a similar way we define
the corresponding objects for $N$, namely $N(\O_y)$ and $N(\K_y)$.

\medskip

We will now recall some useful facts about the structure of
$N(\K_y)$. For the most part, the discussion below applies to
any unipotent algebraic group (in our case it will be $N$).

Note first that $N(\K_y)$ is not only a group ind--scheme (i.e., a
group-like object in the category of ind--schemes), but actually an
ind--group scheme. In other words, $N(\K_y)$ can be represented as a
direct limit of certain group schemes $N^{-k}, k>0$. Moreover, each
$N^{-k}$ is an inverse limit of finite-dimensional unipotent
groups.

Furthermore, for each $k$ we may find a normal subgroup $N^k$ of
$N^{-k}$, such that $N^{k}\subset N(\O_y)$, and $N(\O_y) =
\varprojlim {N(\O_y)/N^{k}}$. We denote by $N_k$
(resp., $N'_k$) the quotient $N^{-k}/N^{k}$ (resp., $N(\O_y)/N^{k}$).

For example, in the case of $GL_2$, when $N={\mathbb G}_a$, $N^k$
(resp., $N^{-k}$) can be chosen to be $t^{-k} \O_y$ (resp.,
$t^k \O_y$), where $t$ is a formal coordinate at $y$.

\sssec{}

So far, our discussion has been local, i.e., we used only the formal
neighborhood of $y$ in $X$. Now the curve $X$ defines a
group ind--subscheme $N_{\out}$ of $N(\K_y)$, such that
$\Hom(S,N_{\out})=\Hom(S\times (X\setminus y), N)$.

Note that any $\GG_a$--bundle over an affine scheme is trivial. Hence,
any $N$--bundle over $(X \backslash y) \times S$, where $S$ is any
$\Fqb$--scheme, can be trivialized locally on $S$ in Zariski
topology. This allows us to obtain all $N$--bundles by ``gluing''
together trivial bundles over $X \backslash y$ and $\D_y$ by a
``transition function'' on $\D_y^\times$, which is an element of
$N(\K_y)$. Thus, informally one can think of $\Bun_N$ as the double
quotient $N_{\out} \backslash N(\K_y)/N(\O_y)$. To avoid quotiening
out by the ind--group $N_{\out}$, we replace $N(\K_y)$ by a large
enough subgroup $N^{-k}$. Then, for large enough $k$, we obtain a
description of $\Bun_N$ as the double quotient
$N_{\on{out},k}\backslash N^{-k}/N(\O_y)$, where
$N_{\on{out},k}=N_{\on{out}}\cap N^{-k}$.

\sssec{} The analogous description of $\BunNft$ for any $T$--bundle
$\F_T$ over $X$ is obtained as follows. Let us choose an embedding
$T\overset{f}\hookrightarrow B$ and consider the induced $B$--bundle
$\F^f_B:=\F_T\overset{T}\times B$. Let $\N^{\F_T}$ be the group scheme
of endomorphisms of $\F^f_B$ that preserve the identification
$\F^f_B/N\simeq \F_T$. 

We have the ind--group scheme $N^{\F_T}(\K_y) =
\Gamma(\D^\times_y,\N^{\F_T})$ and its subgroups $N^{\F_T}(\O_y)$ $=
\Gamma(\D_y,\N^{\F_T})$ and $N^{\F_T}_{\out} = \Gamma(X\setminus y,
\N^{\F_T})$.

Let us fix a trivialization of the restriction of a $T$--bundle
$\F_T$ to $\D_y$:
$$\epsilon:\F^0_T|_{\D_y} \to \F_T|_{\D_y}.$$
Then we can identify $N^{\F_T}(\K_y)$ and $N^{\F_T}(\O_y)$ with
$N(\K_y)$ and $N(\O_y)$, respectively, and consider $N^{\F_T}_{\out}$
as a subgroup of $N(\K_y)$.

Set
$$N^{\F_T}_{\out,k} = N^{\F_T}_{\out} \cap N^{-k} \subset N(\K_y).$$
This is a finite-dimensional unipotent group scheme for any $k$ and we
have a map $N^{\F_T}_{\out,k} \to N_k = N^{-k}/N^k$, which is injective
when $k$ is large enough. Moreover, we have

\sssec{\bf Lemma.}\label{transitivity} {\em For any $k$, a normal
subgroup of finite codimension in $N^k$ acts freely on the quotient
$N^{\F_T}_{\out,k}\backslash N^{-k}$, so the double quotient
$N^{\F_T}_{\out,k}\backslash N^{-k}/N(\O_y)$ is a well-defined
algebraic stack (of finite type). Moreover, we have a natural map
$$N^{\F_T}_{\out,k}\backslash N^{-k}/N(\O_y)\to\BunNft,$$ which is an
isomorphism for $k$ large enough.}

\ssec{The canonical $N(\O_y)$-torsor.}  \label{cantorsor}

\sssec{} \label{exteact}
Let $\F_G$ be an $S$--family of $G$--bundles on $X$.
By making a restriction from $X$ to $\D_y$, we obtain an $S$--family
of $G$--bundles on $\D_y$ (cf. \secref{notation}).

Note that an $S$--family of $G$--bundles on $\D_y$ is the same as an
$G(\O_y)$--bundle over $S$. This construction yields
a canonical $G(\O_y)$--torsor
over $\Bun_G$, which we will denote by $\G_y$.

Analogously, we define the canonical $B(\O_y)$--bundle $\canB_y$ over
$\Bun_B$.  The datum of $\epsilon$ gives us a reduction of
$\canB_y|_{\BunNft}$ to $N(\O_y)$. We denote the corresponding
$N(\O_y)$--bundle over $\BunNft$ by $\canN_y$. The important fact is that
{\it the action of $N(\O_y)$ on $\canN_y$ extends naturally to an
$N(\K_y)$--action.}  Loosely speaking, $N(\K_y)$ acts by changing
the transition function of an $N$--bundle on $\D^\times_y$. A precise
construction of this action is given in the proof of \lemref{Naction}
below.

\sssec{}
Informally, one can say that $\canN_y \simeq N^{\F_T}_{\out} \backslash
N(\K_y)$.  More precisely, \lemref{transitivity} implies the following
result:

\bigskip

\noindent{\bf Corollary.}
{\em For every $k$ we have a natural map $N^{\F_T}_{\out,k}\backslash
N^{-k}\to \canN_y$, which is an isomorphism for $k$ large enough.}

\bigskip

Clearly, the $N^{-k}$ action on $N^{\F_T}_{\out,k}\backslash N^{-k}$
coincides with the one that comes from the above mentioned
$N(\K_y)$--action on $\canN_y$ and the embedding
$N^{-k}\hookrightarrow N(\K_y)$.

\sssec{}
Consider the $N'_k$--bundle $^k\canN_y = \canN_y/N^k$ over
$\BunNft$. \lemref{transitivity} implies:

\sssec{\bf Corollary.}   \label{tower}
{\em For $k$ large enough, we have an isomorphism of $N'_k$--torsors over
$\BunNft$:
$$N^{\F_T}_{\out,k}\backslash N_k\to {}^k\canN_y.$$
In particular, for large $k$, $^k\canN_y$ is an affine scheme,
which is isomorphic to a tower of affine spaces.}

\sssec{Example.}  Let us illustrate the above discussion on the
example of $G=GL(2)$. In this case, we have an isomorphism $\BunNft
\simeq \E_{\L}$ for an appropriate line bundle $\L$
(see \secref{example}). The trivialization $\epsilon$ identifies
$N(\O_y)$ with $H^0(\D_y,\L)$. The stack $\BunNft$ is isomorphic to
$$
H^1(X,\L)/H^0(X,\L) \simeq H^0(X\setminus y,\L) \backslash \K_y/\O_y.
$$
The scheme $\canN_y$ can be
identified with
$$H^1_c(X\setminus y,\L):=\varprojlim
H^1(X,\L(-k\cdot y)) \simeq H^0(X\setminus y,\L) \backslash \K_y.$$

In addition, we have:
$$^k\canN_y \simeq H^1(X,\L(-k\cdot y)) \simeq H^0(X\setminus y,\L)
\backslash \K_y/H^0(\D_y,\L(-k\cdot y)).$$

\sssec{} \label{extension} Now we would like to extend the
$N(\O_y)$--torsor $\canN_y$ from $\BunNft$ to the compactification
$\BunNftb$.  Unfortunately, this does not seem to be possible. The
problem is that those points of $\BunNftb$, for which at least one of
the maps
$$\map^{\lambdach}: \L_{\F_T}^{\lambdach}\to \V^{\lambdach}_{\F_G}$$
has a zero at $y$, do not give rise to genuine $B$--bundles on $\D_y$.

However, $\canN_y$ can be
extended to the open substack $_{y,0}\wt{\Bun}_N^{\ft} \subset
\BunNftb$, since by the definition of $_{y,0}\wt{\Bun}_N^{\ft}$, the data of
$\epsilon$ gives rise to bundle maps
$$\OO|_{\D_y}\to \V^{\lambdach}_{\F_G}|_{\D_y}, \quad \quad
\lambdach\in\check\Lambda^{++},$$ which satisfy the Pl\"ucker
relations. Therefore we obtain a reduction of the
$G(\O_y)$--torsor $\G_y|_{_{y,0}\wt{\Bun}_N^{\ft}}$ to $N(\O_y)$.

We denote the resulting $N(\O_y)$--bundle over
$_{y,0}\wt{\Bun}_N^{\ft}$ by $\CanN_y$.  More precisely, $\CanN_y$
classifies triples $(\F_G,\map,\varphi)$, where $(\F_G,\map) \in
{}_{y,0}\wt{\Bun}_N^{\ft}$ and hence gives rise to a well-defined
$N$--bundle on $\D_y$, and $\varphi$ is a trivialization of this bundle.

The following fact will play an important role in the proof of
\thmmref{one}.

\sssec{\bf Lemma.}    \label{Naction}
{\em The $N(\O_y)$--action on $\CanN_y$ extends naturally to an
$N(\K_y)$--action.}

\sssec{Proof.}

As was mentioned above, the action will come from ``changing
the transition functions over $\D_y^\times$''.

For a point\footnote{By a point we mean, as usual, an $S$--point
for an arbitrary $\Fqb$--scheme $S$.} $(\F_G,\map,\varphi)$ of $\CanN$
and a point $n\in\N(\K_y)$, we will produce a new point
$(\F'_G,\map',\varphi')$ of
$\CanN$ as follows.
We set $(\F'_G,\map')|_{X\setminus y} =
(\F_G,\map)|_{X\setminus y}$. To define the triple
$(\F'_G,\map',\varphi')$ at $y$, we fix an affine neighborhood $X_0$ of
$y$. The trivialization $\varphi$ attaches to $\phi \in
H^0(X_0\setminus y,\V^{\check\la}_{\F_G})$ its Laurent expansion at
$y$, which we denote by $\hat{\phi} \in \V^{\check\la}\otimes
\K_y$. We define $\F'_G$ on $X_0$ by the following condition: the sections
of $\V^{\check\la}_{\F'_G}$ on $X_0$ are precisely
the sections $\phi \in
H^0(X_0\setminus y,\V^{\check\la}_{\F_G})$ such that $n \cdot \hat{\phi}
\in \V^{\check\la}\otimes \O_y$. It remains to check that the  map
$\map'$ indeed maps
$\L_{\F_T}^{\check\la}$ to $\V^{\check\la}_{\F'_G}$. This follows from the
fact that the image of  the map $\map^{\check\la}|_{\D^\times_y}:
\L_{\F_T}^{\check\la}|_{\D^\times_y}
\to\V^{\check\la}_{\F_G}|_{\D^\times_y}$ is
$N(\K_y)$--invariant.

Finally, note that the  $G$--bundle
$\F'_G$ comes with a tivialtization over $\D_y$ and by
construction this trivialization is compatible with $\map'$. This gives
the trivialization $\varphi'$.

\bigskip

Note that the $N(\K_y)$--action does not preserve the projection
$\CanN_y \to\BunNft$.

\ssec{Affinness.}

In this section we prove the following statement, which will be used
in the proof of \propref{stability}.

\sssec{\bf Proposition.}  \label{affine} {\em The open embedding
$_{\xl,\nul}\BunNft \hookrightarrow _{\xl,\nul}\BunNftb$ is affine.}

\sssec{}

Recall that we have an isomorphism \eqref{natiso}, which identifies
$_{\xl,\nul}\BunNft$ and $\Bun_N^{\F'_T}$ as well as their
closures. Therefore to prove \propref{affine} it suffices to show that
the embedding $\BunNft \hookrightarrow \BunNftb$ is affine.

For $\mu \in \Lambda^+$, introduce an open substack
$\BunNftbm$ of $\BunNftb$. It classifies those pairs
$(\F_G,\map)$ for which the torsion part of the quotient sheaf
$\V^{\lambdach}_{\F_G}/\on{Im}(\map^{\lambdach})$ has length less than
or equal to $\langle \mu,\lambdach\rangle$, for all
$\lambdach\in\check\Lambda^{++}$.

Clearly, $\ol{\Bun}_{N,0}^{\ft}=\BunNft$ and
$\BunNftbm \subset \ol{\Bun}_{N,\nu}^{\ft}$ if
$\nu-\mu\in\Lambda^+$. For $y\in X$ let $_{y,0}\wt{\Bun}^{\F_T}_{N,\mu} ={}
_{y,0}\wt{\Bun}^{\F_T}_{N} \cap \BunNftbm$.

It is clear that that every $\Fqb$-point of $\BunNftb$ belongs to
$_{y,0}\wt{\Bun}^{\F_T}_{N,\mu}$ for some $y$ and $\mu$. Hence, it is
sufficient to show that the embedding of $\BunNft$ into every
$_{y,0}\BunNftm$ is affine.

\bigskip

Now \propref{affine} follows by combining \corref{tower} and the
following result. Denote the $N'_k$--bundle $\CanN_y/N_k$ by
$^k\CanN_y$.

\sssec{\bf Proposition.}  \label{scheme} {\em For fixed $\F_T$, $\mu$
and $y$, the restriction of $^k\CanN_y$ to
$_{y,0}\wt{\Bun}^{\F_T}_{N,\mu}$ is a scheme when $k$ is large
enough.}

\sssec{Proof.}

First, we claim that the stack $_{y,0}\wt{\Bun}^{\F_T}_{N,\mu}$ is of
finite type. This follows from:

\noindent {\bf Lemma.}  {\em The image of $_{y,0}\BunNftm$ under the
projection $\p:\BunNftb\to\Bun_G$ is contained in an open substack $U$
of $\Bun_G$ of finite type.}

\bigskip

\noindent{\em Proof of the lemma.}  Let $\check\la\in\check\Lambda^{++}$ be a
regular weight and let $\{\check\eta\}$ be the collection of weights
of $\V^{\check\la}$. We can choose $\nu \in \check\Lambda$ such that
$$\langle\nu,\check\la\rangle > \underset{\check{\eta}; 0\leq \mu'\leq
\mu}{\on{max}} \langle \mu'+\deg(\F_T),\check\eta\rangle.$$

Let $\F_G$ be a $G$--bundle in the image of $_{y,0}\BunNftm$ under the
projection $\p:\BunNft\to\Bun_G$. By definition, it admits a
$B$--structure of degree $\deg(\F_T)+\mu'$ for $0\leq\mu'\leq \mu$
(here by degree of a $B$--bundle we understand the degree of the
corresponding $T$--bundle, which is an element of $\Lambda$).  Then
the associated bundle $\V^{\check\la}_{\F_G}$ admits a filtration,
whose successive quotients are line bundles of degrees $\langle
\mu'+\deg(\F_T),\check\eta\rangle$.

Therefore $\V^{\check\la}_{\F_G}$ does not admit a line subbundle of
degree $\langle \nu',\check\la\rangle$ for any $\nu'\geq\nu$. Hence,
$\F_G$ does not admit a $B$--structure of degree greater than or equal
to $\nu$. But it is well-known that the open substack of $\Bun_G$,
which classifies such $G$--bundles, is of finite type. \qed

Recall the canonical $G(\O_y)$--bundle $\G_y$ over $\Bun_G$ from
\secref{exteact}. For $i>0$, denote by $G^{i}$ the $i$--th congruence
subgroup of $G(\O_y)$ and by $G_i$ the quotient $G(\O_y)/G^i$.  The
stack--theoretic quotient $^i\G_y:=\G_y/G^{i}$ is a
principal $G_i$--bundle over $\Bun_G$. It is well-known that for any
open substack $U\subset\Bun_G$ of finite type, there exists an integer
$i>0$ such that the restriction of $^i\G_y$ to $U$ is a scheme.

Since $\p$ is representable, we obtain that on the one hand,
$^i\G_y|_{_{y,0}\wt{\Bun}^{\F_T}_{N,\mu}}$ is a scheme. On the other
hand, by definition, the induced $G(\O_y)$--bundle $\CanN_y
\underset{N(\O_y)}\times G(\O_y)$ can be canonically identified with
$\G_y|_{_{y,0}\wt{\Bun}^{\F_T}_{N,\mu}}$.

Now, for every $i>0$ we can find an integer $k>0$ such that
$N^k\subset N(\O_y)\cap G^i$. Therefore for a given $i$ and large
enough $k$, we have an isomorphism of $G_i$--bundles:
$$^k\CanN_y\underset{N'_k}\times
G_i\simeq {}^i\G_y|_{_{y,0}\BunNftm}.$$
Hence for $i$ and $k$ large enough, $^k\CanN_y\underset{N'_k}\times
G_i$ is a scheme. Since the stack $^k\CanN_y$
admits a representable morphism to $^k\CanN_y\underset{N'_k}\times
G_i$, it is also a scheme for $i$ and $k$ large enough.

This completes the proof of \propref{scheme} and therefore of
\propref{affine}.

\section{Sheaves}   \label{category}

\ssec{Evaluation morphisms and perverse sheaves}    \label{evaluation}

\sssec{}  \label{Brep}

For a fixed $\F_T\in\Bun_T$ consider the affine space: $\prod_{\i\in
\I} H^1(X,\L_{\F_T}^{\alphach_\i})$.

We construct a natural morphism:
$$\on{ev}:\BunNft\to \prod_{\i\in \I}
H^1(X,\L_{\F_T}^{\alphach_\i}).$$

First, let us consider the case $G=GL_2$.  Let $\L_1$ and $\L_2$ be
as in \secref{example}. Then
$\L_{\F_T}^{\alphach}\simeq\L:=\L_1\otimes\L_2^{-1}$. The map
$\on{ev}$ is just the composition
$$\BunNft\simeq\E_{\L}\to H^1(X,\L).$$

To define the morphism $\on{ev}$ for general $G$ we construct for
every $\i\in\I$ a two-dimensional representation $\V^\i$ of the group
$B$. Its restriction to $T \subset B$ is the direct sum $\one \oplus
\one^{\alphach_{\i}}$ of the trivial one-dimensional representation
$\one$ and the one--dimensional representation corresponding to the
character $B\to T\overset{\alphach_\i}\to {\mathbb G}_m$. The subgroup
$N \subset B$ maps $\one$ to $\one^{\alphach_{\i}}$ via the
projection
$$N\twoheadrightarrow N/[N,N]\to \on{Span}\{\alphach_\i\}.$$ Thus, we
have a short exact sequence of $B$--modules
$$0\to \one^{\alphach_\i}\to \V^\i\to\one\to 0.$$

By definition, a point $(\F_G,\map)\in\BunNft$ defines a $B$-bundle
$\F_B$ on $X$, such that $\V^\i_{\F_B}$ fits into a short exact
sequence:
$$0\to \L_{\F_T}^{\alphach_\i}\to \V^\i_{\F_B}\to {\mc O}_X \to 0.$$

Hence we obtain a morphism of stacks $\BunNft\to
\E_{\L_{\F_T}^{\alphach_\i}}$ (the latter is defined in
\secref{example}). By composing it with the canonical map
$\E_{\L_{\F_T}^{\alphach_\i}}\to H^1(X,\L_{\F_T}^{\alphach_\i})$ we
obtain a morphism
$$\on{ev}_\i:\BunNft\to H^1(X,\L_{\F_T}^{\alphach_\i}).$$
Finally, we set $\on{ev}=\prod_{\i\in\I} \on{ev}_\i$.

\sssec{}  \label{charactersheaves}

To define our sheaves on $\BunNft$ we need to choose additional data.

Suppose that for every $\i\in\I$ we are given an embedding
$$\varpi_\i:\L_{\F_T}^{\alphach_\i}\hookrightarrow\Omega.$$
Let $\varpi$ denote the collection $\{\varpi_\i\},\,\forall\,\i\in\I$.

To such $\varpi$ we assign its conductor, $\on{cond}(\varpi)$, which
is a divisor on $X$ with values in the semi-group
$\Lambda_{G_{\on{ad}}}^{++}$ of dominant co-weights of the group $G_{\on{ad}}$,
where $G_{\on{ad}}:=G/Z(G)$. Namely, we set $\langle
\on{cond}(\varpi),\alphach_\i\rangle$ equal to the divisor of zeroes
of the map $\varpi_\i$ for all $\i\in\I$.

For a point $x\in X$, the conductor of $\varpi$ at $x$, denoted by
$\on{cond}_x(\varpi)$, is by definition an element of
$\Lambda_{G_{\on{ad}}}^{++}$ equal to the value of $\on{cond}(\varpi)$ at
$x$.

\sssec{}    \label{ArtS}

Let us fix the data of $\varpi$.  We define a morphism
$\on{ev}_\varpi:\BunNft\to \GG_a$  as the composition
$$\BunNft\overset{\on{ev}_\i} \longrightarrow \prod_{\i\in\I}
H^1(X,\L^{\alphach_\i}) \overset{\varpi}\longrightarrow
\prod_{\i\in\I} H^1(X,\Omega)\simeq
\GG_a^\I\overset{\on{sum}}\longrightarrow \GG_a,$$ where the last arrow
is the map $(a_1,...,a_r)\to a_1+...+a_r$.

Write $\J_\psi$ for the Artin-Schreier sheaf on
$\GG_a$ corresponding to a fixed additive character $\psi:\Fq\to\Ql^*$
and let
$\Psi_{\varpi}$ be the local system on $\BunNft$ defined by
$$\Psi_{\varpi} := \on{ev}_{\varpi}^*(\J_\psi)[d_N],$$
where $d_N = \on{dim} \BunNft$.

Let $\Psib_{\varpi}$ be the perverse sheaf on $\BunNftb$, which is the
Goresky-MacPherson extension of $\Psi_{\varpi}$:
$$\Psib_{\varpi}:=j_{!*}(\Psi_{\varpi}).$$

We also define {\it complexes} of perverse sheaves
$$
\Psi_{\varpi}{}_{!} := j_{!}(\Psi_{\varpi}), \quad \quad
\Psi_{\varpi}{}_{*} := j_{*}(\Psi_{\varpi})
$$
on $\BunNftb$.

\begin{thmm}    \label{oneprime}

The canonical maps
$$\Psi_{\varpi}{}_{!}\to \Psib_{\varpi} \to \Psi_{\varpi}{}_{*}$$ are
isomorphisms.

\end{thmm}

In other words,  the $*$--restriction of
$\Psib_{\varpi}$ to the complement of $\BunNft$ in $\BunNftb$
vanishes.

\ssec{The definition of the category of sheaves}

\sssec{}

Next we reformulate \thmmref{oneprime} using the ind-stack
$_{\xl,\infty}\BunNftb$.

Let $(\xl,\nul)$ be as in \secref{variants} and consider the stack
$_{\xl,\nul}{\Bun}{}^{\F_T}_N$. Denote by $\on{pr}$ the natural
projection $\Lambda\to \Lambda_{G_{\on{ad}}}$. Let us assume that for
every $k$, the sum
$$\on{cond}_{x_k}(\varpi)+\on{pr}(\nu_k)$$ is still a
dominant coweight of $G_{\on{ad}}$.

We have an identification of stacks
$$_{\xl,\nul}{\Bun}{}^{\F_T}_N\simeq \Bun^{\F'_T}_N,$$
where $\F'_T:=\F_T\left( -\sum_k \nu_k\cdot x_k \right)$. The above
condition on $\varpi$ ensures that the corresponding meromorphic maps
$$\varpi'_\i:\L^{\alphach_\i}_{\F'_T}\to \Omega$$
are regular. Hence we obtain a morphism $\on{ev}_{\varpi'}:
\Bun^{\F'_T}_N \arr \GG_a$, and hence a morphism
$$\on{ev}^{\xl,\nul}_{\varpi}:
{}_{\xl,\nul}{\Bun}{}^{\F_T}_N\to\GG_a.$$

Let $\Psi^{\xl,\nul}_{\varpi} = (\on{ev}^{\xl,\nul}_{\varpi})^*({\mc
I}_\psi)[d_N]$ be the corresponding local system on
$_{\xl,\nul}{\Bun}{}^{\F_T}_N \subset{} _{\xl,\infty}\BunNft$. Let us
denote by $\Psib^{\xl,\nul}_{\varpi}$ be the Goresky-MacPherson
extension of $\Psi^{\xl,\nul}_{\varpi}$, and
$\Psi^{\xl,\nul}_{\varpi}{}_{!}$ and $\Psi^{\xl,\nul}_{\varpi}{}_{*}$
be the corresponding complexes of perverse sheaves on
$_{\xl,\infty}\BunNftb$.

\sssec{\bf Definition.}  Let us fix data of $\varpi$, such that
$\on{cond}(\varpi)=0$ (note that this means that we fix  isomorphisms
$\L^{\alphach_{\i}}_{\F_T} \simeq \Omega$). We define the {\em
Whittaker category} $\W^{\xl}_\varpi$ as the full abelian subcategory
of the category of all perverse sheaves on $_{\xl,\infty}\BunNftb$,
which consists of the perverse sheaves whose irreducible subquotients
are the sheaves $\Psib^{\xl,\nul}_\varpi$, with $\nul$ such that
$\nu_k\in\Lambda^{++}$ for all $k=1,...,n$.

\bigskip

The following is the main result of this paper.

\begin{thmm}    \label{one}

\begin{enumerate}

\item[(1)] The category $\W^{\xl}_\varpi$ is semi-simple.

\item[(2)] The canonical maps $\Psi^{\xl,\nul}_\varpi{}_{!}\to
\Psib^{\xl,\nul}_\varpi\to \Psi^{\xl,\nul}_\varpi{}_{*}$ are
isomorphisms.

\end{enumerate}

\end{thmm}

\thmmref{one} will be proved in \secref{proofs}. As was already
mentioned in \secref{intr}, \thmref{oneprime}
will be be obtained as a rather formal consequence of the explicit
computation of the action of the Hecke functors on $\W_\varpi^{\xl}$,
which will be introduced in the next section.

\sssec{} \label{independence}
We claim that \thmmref{one}(2) and \thmmref{oneprime} are equivalent.
Clearly, \thmmref{one}(2) is a particular case of
\thmmref{oneprime}. Let us explain how to derive \thmmref{oneprime} from
\thmmref{one}(2).

Suppose first that the center $Z(G)$ of $G$ is connected. This means
that the map $\on{pr}:\Lambda\to \Lambda_{G_{\on{ad}}}$ is surjective. In
this case, starting with an arbitrary $\varpi$, there exists a
$\Lambda$--valued divisor $D=\sum_k \mu_k \cdot x_k$ on $X$, such that
$\on{pr}(D)=\on{cond}(\varpi)$.

Thus, if we set $\F'_T:=\F_T(D)$, we will have $\on{cond}(\varpi')=0$
for the corresponding maps $\varpi'_\i: \L^{\alphach_\i}_{\F'_T}\to
\Omega$. Moreover, we can then identify $_{\xl,\nul}{\Bun}{}^{\F_T}_N$
with $_{\xl,\nul+\mul}{\Bun}^{\F'_T}_N$, and identify their closures
in $_{\xl,\infty}\BunNftb$ and $_{\xl,\infty}\ol{\Bun}^{\F'_T}_N$,
respectively. The morphisms $\on{ev}_{\varpi}:
_{\xl,\nul}{\Bun}^{\F_T}_N \to \GG_a$ and $\on{ev}_{\varpi'}:
_{\xl,\nul+\mul}{\Bun}{}^{\F'_T}_N \to \GG_a$ coincide under this
identification, as do the corresponding perverse sheaves.

Therefore we see that in the case when $Z(G)$ is connected
\thmmref{oneprime} and \thmmref{one}(2) are equivalent.

\sssec{} Let now $G$ be an arbitrary reductive group, such that
$[G,G]$ is simply-connected, and consider the group
$G_1:=G\underset{Z(G)}\times T$. Evidently, $[G_1,G_1]$ is also
simply-connected. In addition, $G_1$ has a connected center:
$Z(G_1)=T$ and \thmmref{oneprime} for $G_1$ follows from \thmmref{one}.

Let $T_1$ denote the Cartan subgroup of $G_1$; we have an embedding
$T\hookrightarrow T_1$ and let $\F_{T_1}$ denote the induced
$T_1$--bundle over $X$.

We have a natural map of stacks $\BunNftb\to
\overline{\Bun}^{\F_{T_1}}_{N_1}$, which is easily seen to be an
isomorphism.

Moreover, the data of $\varpi$ for $G$ defines the corresponding data
$\varpi_1$ for $G_1$ such that the sheaves
$\Psi_{\varpi}{}_{!}$, $\Psib_{\varpi}$ and $\Psi_{\varpi}{}_{*}$
match with the corresponding perverse sheaves for $G_1$. Hence,
\thmmref{oneprime} for $G$ will follow from \thmmref{one}(2) for $G_1$.

\section{Hecke functors} \label{hecke}

\ssec{The affine Grassmannian}

\sssec{} \label{affgrass}

Let $x\in X$ be an $\Fqb$--point, $\D_x$ the formal disc around $x$ in
$X$, and $t$ a formal coordinate at $x$.  The choice of the
formal coordinate allows us to identify $\O_x$ with $\Fqb[[t]]$ and
$\K_x$ with $\Fqb((t))$.

Recall the group scheme $G(\O_x)$ and the ind--group scheme
$G(\K_x)$ that classify maps $\D_x\to G$ and $\D_x^\times\to G$,
respectively.

The affine Grassmannian $\Gr$ is an ind-scheme defined as the quotient
$G(\K_x)/G(\O_x)$. In other words, for a scheme $S$, $\Hom(S,\Gr)$ is
the set of pairs $(\F_G,\beta)$, where $\F_G$ is an $S$--family of
$G$--bundles on $\D_x$ and $\beta$ is a trivialization of the
corresponding family of bundles on $\D_x^\times$.

Alternatively, one can define $\Gr$ using the global curve $X$:
$\Hom(S,\Gr)$ is a set of pairs $(\F_G,\beta)$, where $\F_G$ is a
$G$--bundle over $X\times S$ and $\beta$ is a trivialization of the
restriction of $\F_G$ to $(X\setminus x)\times S$.

Evidently, to a data $(\F_G,\beta)$ ``over $X$'' one can attach a data
$(\F_G,\beta)$ ``over $\D_x$'' and the fact that the two definitions
coincide is a theorem due to Beauville and Laszlo \cite{BL}.

\sssec{} The first description of $\Gr$ implies that it carries a
natural action of $G(\K_x)$ and, in particular, of $G(\O_x)$. For
$\la\in\PL^+$ let $\Gr^\la$ be the $G(\O_x)$-orbit of the point
$\la(t) \cdot G(\O_x)\in \Gr$, where $\la(t)\in T(\K_x) \subset
G(\K_x)$. It is easy to see that $\Gr^\la$ is independent of the
choice of $T\subset G$ and of the uniformizer $t\in\O_x$. Furthermore,
$\on{dim}({\Gr}^\la)=\langle\la,2\rhoc\rangle$ and
 $\Gr^\mu \subset \ol{\Gr}^\la$ if and only of
$\lambda\geq\mu$, i.e.,  when $\la-\mu$ is a sum of positive roots of
$\GL$ -- see, for example, \cite{Lu2}. The schemes
$\Gr^\la$ form a stratification of the underlying reduced scheme of
$\Gr$. 

\ssec{The category $\on{P}_{G(\O_x)}(\Gr)$} \label{defcat}

\sssec{}

Since the closures of $G(\O_x)$--orbits on $\Gr$ are
finite-dimensional, the abelian category $\on{P}_{G(\O_x)}(\Gr)$ of
$G(\O_x)$--equivariant perverse sheaves on $\Gr$ is well-defined.  By
definition, every object of this category is  supported over a finite
union of schemes of the form $\ol{\Gr}^\la$.

It is well-known that for every $\Fqb$-point of $\Gr$, its stabilizer
in $G(\O_x)$ is connected.  Therefore, the irreducible objects in
$\on{P}_{G(\O_x)}(\Gr)$ are the Goresky--MacPherson extensions of
constant local systems (appropriately shifted) on the strata $\Gr^\la$; we
will denote these sheaves by ${\mathcal A}_\la$.

\sssec{\bf Theorem.}(see \cite{Lu2}) \label{semisimplicity}
{\em The category $\on{P}_{G(\O_x)}(\Gr)$ is semi-simple.}

\bigskip

We remark that the argument given in \secref{proofofsemisimplicity}
yields an alternative proof of this theorem.

\sssec{}  \label{intrconv}

Next we will describe the convolution operation
$\on{P}_{G(\O_x)}(\Gr)$. This construction is nothing but a
sheaf-theoretic analogue of the definition of the algebra structure on
$\HH$.

\medskip

Consider the ind-scheme $\Conv:=G(\K_x)\underset{G(\O_x)}\times
G(\K_x)/G(\O_x)$. It is easy to see that this scheme represents the
following functor: $\Hom(S,\Conv)$ is the set of quadruples
$(\F^1_G,\beta^1,\F_G,\beta)$, where $\F^1_G$ and $\F_G$ are
$G$--torsors over $\D_x$, $\beta^1$ is an identification
$\F^1_G|_{\D_x^\times}\to
\F_G|_{\D_x^\times}$ and $\beta$ is a trivialization
$\F_G|_{\D_x^\times}\to
\F^0_G|_{\D_x^\times}$.

There is an obvious projection $p_1:\Conv\to\Gr$ which sends the
quadruple \newline
$(\F^1_G,\beta^1,\F_G,\beta)$ to $(\F_G,\beta)$; it
is the projection on the first factor:
$$G(\K_x)\underset{G(\O_x)}\times G(\K_x)/G(\O_x) \to G(\K_x)/G(\O_x)
= \Gr.$$ Moreover, this projection realizes $\Conv$ as a fibration
over $\Gr$ with a typical fiber isomorphic again to $\Gr$.

More precisely, $\Conv$ is a bundle associated to the
$G(\O_x)$--scheme $\Gr$ and $G(\O_x)$--torsor $G(\K_x)\to
G(\K_x)/G(\O_x)=\Gr$.

Thus, starting from a perverse sheaf $\S_1$ on $\Gr$ and a
$G(\O_x)$--equivariant perverse sheaf $\S_2$ on $\Gr$ we can construct
their twisted external product $\S_1\widetilde\boxtimes\S_2$, which is
a perverse sheaf on $\Conv$ (see \secref{notation}).

Note that the product map gives rise to a second projection
$$\Conv=G(\K_x)\underset{G(\O_x)}\times G(\K_x)/G(\O_x)\to
G(\K_x)/G(\O_x)=\Gr.$$ On the level of the corresponding functors,
this map sends a quadruple $(\F^1_G,\beta^1,\F_G,\beta)$ as above to
$(\F_G^1,\beta')$, where $\beta'$ is obtained as a composition
$\beta\circ \beta^1$. We will denote this projection by $p_2$.

For two perverse sheaves $\S_1$ and $S_2$ on $\G$ with $\S_2$ being
$G(\O_x)$--equivariant, we denote by $\S_1\star\S_2$ the complex
$p_2{}_{!}(\S_1\widetilde\boxtimes\S_2)$. By construction, if
$\S_1$ is also $G(\O_x)$--equivariant, then so is $\S_1\star\S_2$.

\sssec{\bf Theorem.}(\cite{Lu2}) {\em For $\S_1,\S_2\in
\on{P}_{G(\O_x)}(\Gr)$, the complex $\S_1\star\S_2$ is a perverse
sheaf.}

\bigskip

Actually, even more is true: $\S_1\star\S_2$ is a perverse sheaf for
$\S_1 \in \on{P}_{G(\O_x)}(\Gr)$ and any perverse sheaf
$\S_2$. However, we will not need this stronger statement.

\sssec{} Thus, $\S_1,\S_2\mapsto \S_1\star\S_2$ defines a binary operation
on $\on{P}_{G(\O_x)}(\Gr)$. It is easy to see from the definition of
$\star$ that it admits a natural associativity constraint, which makes
$\on{P}_{G(\O_x)}(\Gr)$ into a monoidal category.

In addition, the inversion $g\mapsto g^{-1}$ on $G(\K_x)$ induces a
covariant self--functor on $\on{P}_{G(\O_x)}(\Gr)$, which we will
denote by $\S\mapsto \ast\S$.

The following theorem can be viewed as a sheaf--theoretic version of
the Satake isomorphism (see \thmref{Satake}).

\sssec{\bf Theorem.}  \label{sheaf-Satake} {\em The category
$\on{P}_{G(\O_x)}(\Gr)$ admits a commutativity constraint, i.e.,
it is a tensor category. Moreover, as such, it is equivalent
to the category ${{\mathcal R}ep}(\GL)$ of finite-dimensional
representations of $\GL$ in such a way that}

\begin{enumerate}

\item[(1)] {\em The object $\A_\la\in\on{P}_{G(\O_x)}(\Gr)$ goes to
$V^\la$.}

\item[(2)] {\em Let $\DD$ be the Verdier duality functor on
$\on{P}_{G(\O_x)}(\Gr)$. The functor $\S\to \DD(\ast\S)$ corresponds to
the contragredient duality functor.}

\end{enumerate}

\sssec{Remark.}
To be precise, this result has been proved in \cite{Gi,MV} over the
ground field $\CC$ (in this setting, this isomorphism was conjectured
by V.~Drinfeld; see also \cite{Lu2}). But the proof outlined in
\cite{MV} can be generalized to the case of the ground field $\Fqb$;
see also \cite{Ngo1}.

\ssec{The Hecke action on $_{x,\infty}\BunNftb$}

\sssec{\bf Definition.}
For a point $x \in X$ the {\em Hecke correspondence} stack $\H_x$
is defined as follows. For an $\Fqb$--scheme $S$, $\on{Hom}(S,\H_x)$ is a
groupoid, whose objects are triples $(\F_G,\F'_G,\beta)$, where $\F_G$
and $\F'_G$ are $G$--bundles over $X \times S$ and $\beta$ is an
isomorphism
$$\beta: \F_G|_{(X\setminus x) \times S} \simeq \F'_G|_{(X\setminus x)
\times S}.$$
The definition of morphisms in $\on{Hom}(S,\H_x)$ is straightforward.

Let $\hl$ and $\hr$ be the projections $\H_x \arr \Bun_G$, sending
$(\F_G,\F'_G,\beta)$ to $\F_G$ and $\F'_G$, respectively.

\sssec{} Recall the $G(\O_x)$--bundle $\G_x$ over $\Bun_G$. It
follows from the definitions that each of the projections $\hl$ and
$\hr$ realizes $\H_x$ as a fibration over $\Bun_G$ attached to
$\G_x$ and the $G(\O_x)$--scheme $\Gr$:
\begin{equation}    \label{hx}
\H_x \simeq \G_x \underset{G(\O_x)}{\times} \Gr.
\end{equation}

Thus, if $\S$ is an object of $\on{P}_{G(\O_x)}(\Gr)$, we can use
either $\hr$ or $\hl$ to produce objects $\S^l$ and $\S^r$ in $\Sh(\H_x)$
by taking the twisted external product of $\S$ with the constant sheaf
on $\Bun_G$.  Note that there is a natural isomorphism
\begin{equation}  \label{left-right}
(\ast\S)^l\simeq\S^r.
\end{equation}

\medskip

Let $\Sh(\Bun_G)$ be the derived category of $\Ql$--sheaves on
$\Bun_G$. It is well-known that $\on{P}_{G(\O_x)}(\Gr)$ ``acts'' on
$\Sh(\Bun_G)$ by convolution (this action is used in the
definition of ``Hecke eigensheaves'' in the geometric Langlands
correspondence). We will need an analogous action of
$\on{P}_{G(\O_x)}(\Gr)$ on the derived category $\Sh(\bnx)$ of
$\Ql$--sheaves on $\bnx$.

\sssec{}

Let $Z$ denote the fiber product
\begin{equation}    \label{cart}
Z = \H_x \underset{\Bun_G}\times{}_{x,\infty} \BunNftb,
\end{equation}
where the morphism $\H_x\to\Bun_G$ that we use in the above formula is
the projection $\hr$.

In other words, $Z$ is a fibration over $\bnx$ that corresponds to the
$G(\O_x)$--torsor $\G_x|_{\bnx}$ and the $G(\O_x)$--scheme $\Gr$.

\sssec{\bf Proposition.}    \label{secondprojection}
{\em There exists a second projection
$Z\overset{'\hl}\longrightarrow{}_{x,\infty}\BunNftb$ which renders
the diagram
\begin{equation}  \label{comdiag}
\CD
_{x,\infty}\BunNftb   @<{'\hl}<<   Z  @>{'\hr}>>   _{x,\infty}\BunNftb  \\
@V{\p}VV      @V{'\p}VV        @V{\p}VV   \\
\Bun_G  @<{\hl}<<  \H_x @>{\hr}>>  \Bun_G
\endCD
\end{equation}
commutative. Moreover, the left square in the above diagram is
Cartesian as well.}

\sssec{Proof.}
By definition, the stack $\bnx$ classifies pairs $(\F_G,\map)$, where
where $\map$ is a collection of maps
$$\map^{\lambdach}:\L_{\F_T}^{\lambdach}\to
\V^{\lambdach}_{\F'_G}(\infty\cdot x)$$
subject to the Pl\"ucker relations.

Therefore, the stack $Z$ classifies triples
$(\F_G,\F'_G,\map,\beta)$,
with $\F_G$ and $\map$ are as above and $\beta$ is an isomorphism
$\F_G|_{X\setminus x}\simeq \F'_G|_{X\setminus x}.$

Hence, from $\map$ and $\beta$ we obtain a
system of maps
$$\map'{}^{\lambdach}:\L_{\F_T}^{\lambdach}\to
\V^{\lambdach}_{\F'_G}(\infty\cdot x),$$
which satisfy the Pl\"ucker relations.

Let us set $'\hl(\F_G,\F'_G,\map,\beta) = (\F_G,\map')$, with the
collection $\map' = \{ \map'{}^{\lambdach} \}$ defined above. It is
clear that the left square of the above diagram is commutative and
Cartesian.

\sssec{}

The morphism ${}'\hr$ realizes $Z$ as a fibration $\G_x
\underset{G(\O_x)} \times \Gr$. Therefore given $\S \in
\on{P}_{G(\O_x)}(\Gr)$ and $\T \in \Sh(\bnx)$, using the notation
introduced in \secref{notation}, we can form the twisted tensor
product $(\T \tboxtimes \S)^r$. Note that because ${}'\hl$ is  proper on
the support of $(\T \tboxtimes \S)^r$, we have
$$
{}'\hl_{*}((\T \tboxtimes \S)^r)\ = \ {}'\hl_{!}((\T
\tboxtimes\S)^r).$$
We then set:
\begin{equation}
\T\star \S:={}'\hl_{!}((\T \tboxtimes \S)^r) \ = \ {}'\hl_{*}((\T
\tboxtimes \S)^r).
\end{equation}

Analogously, considering $Z$ as a fibration $\G_x \underset{G(\O_x)}
\times \Gr$ corresponding to the morphism ${}'\hl$, we can construct
the twisted tensor product $(\T \tboxtimes \S)^l$. Then we obtain a
functor ``from left to right'':
$$\T,\,\S\mapsto \S\star\T = {}'\hr_!((\T\widetilde{\boxtimes}(*\S))^l).$$

\sssec{Properties of the action.}
The Hecke functors defined above
$$\Perv_{G(\O_x)}(\Gr)\times \Sh(\bnx)\to\Sh(\bnx)$$
are compatible with the tensor structure on $\Perv_{G(\O_x)}(\Gr)$
in the following sense:

\sssec{\bf Lemma.}    \label{conv}
{\em For any $\S_1,\S_2\in \Perv_{G(\O_x)}(\Gr)$ we have functorial
isomorphisms
\begin{align*}
& (\T \star \S_1) \star \S_2 \simeq \T \star (\S_1 \star \S_2), \\
&\S_1 \star (\S_2 \star \T)) \simeq (\S_1\star\S_2) \star \T,\\
\end{align*}
such that the pentagon identity holds for three-fold convolution
products.}

\bigskip

Moreover, we also have the following proposition which follows in a
straightforward way from the definitions..

\sssec{\bf Proposition.}    \label{properties}

\begin{enumerate}

\item[(1)] {\em The Hecke functors commute with  Verdier duality in
the sense that $\DD(\S \star \T)) \simeq \DD(\S) \star \DD(\T)$, and
$\DD(\T \star \S)) \simeq \DD(\T) \star \DD(\S)$.}

\item[(2)] {\em For a fixed $\S\in \on{P}_{G(\O_x)}(\Gr)$, the
functors $\T\mapsto \T \star \S$ and $\T\mapsto \DD(\S) \star \T$ from
$\Sh(\bnx)$ to itself are mutually both left and right adjoint.}

\item[(3)] {\em There is a functorial isomorphism $\T\star\S\simeq
(*\S)\star \T$.}

\end{enumerate}

\ssec{Action on the canonical sheaves}

\sssec{} Now we take as $\S$, the perverse sheaf ${\mathcal A}_\la$,
which is the Goresky--MacPherson extension of the constant sheaf
(appropriately shifted) on
$\Gr^\la$. Note that we have: $*\A_\la = \A_{-w_0(\la)}$.

The following statement, which will be proved in \secref{heckecontd},
is a crucial step in our proof of \thmmref{one}(1).

\begin{thmm}    \label{two}
Suppose that $\on{cond}_x(\varpi)=0$. Then
$$\Psib_\varpi^{x,0}\star\A_\la\simeq \Psib^{x,\la}_\varpi.$$
\end{thmm}

\sssec{}

According to \lemref{conv} and \thmref{sheaf-Satake}, for any
$\mu\in\Lambda^{+}$, there exists a functorial isomorphism
$$(\T\star\A_{\lambda}) \star \A_{\mu}  \simeq
\oplus_{\nu\in\Lambda^+} (\T\star\A_\nu)\otimes
\on{Hom}_{\GL}(V^\nu,V^\lambda\otimes V^\mu).$$
Therefore we obtain:

\sssec{\bf Corollary.}    \label{decomp}
{\em Suppose that $\on{cond}_x(\varpi)=0$. Then}
\begin{equation}    \label{decomp1}
\Psib^{x,\mu}_\varpi\star\A_\la \simeq \oplus_{\nu\in\Lambda^{+}}
\Psib^{x,\nu}_\varpi \otimes \on{Hom}_{\GL}(V^\nu, V^\lambda \otimes
V^\mu)
\end{equation}

\sssec{The category $P_{N(\K_x)}^\chi(\Gr)$.}    \label{right category}

Let us explain why \thmmref{one} and \thmmref{two} allow us to construct
a category $P_{N(\K)}^\chi(\Gr)$ which satisfies the properties listed
in \secref{expprop}.

Choose data of $\varpi$ with $\on{cond}(\varpi)=0$ (in fact,
$\on{cond}_x(\varpi)=0$ would suffice). Set
$$
P_{N(\K)}^\chi(\Gr):=\W^x_\varpi.
$$

Then \thmmref{two} implies that the Hecke action of
$\on{P}_{G(\O_x)}(\Gr)$ on $\Sh(\bnx)$
preserves $\W^x_\varpi$ and that $\Psib^{x,0}_\varpi\star \A_\la\simeq
\Psib^{x,\la}_\varpi$.

This gives us the first two properties of $P_{N(\K)}^\chi(\Gr)$. The
third property is insured by \thmmref{one}.

We will see in \secref{conj} that the computation of stalks of
$\Psib^{x,\mu+\nu}_\varpi\star\A_\la$ is equivalent to the
computation of the cohomology \eqref{gencs}.

\bigskip

Our plan now is the following. In \secref{proofs} we derive
\thmmref{one} from \thmmref{two}. Then in \secref{heckecontd} we prove
\thmmref{two}. Finally, we use these results in \secref{conj} to
prove \thmmref{gencs1}.

\section{Proof of \thmmref{one}}  \label{proofs}

\ssec{Proof of \thmmref{one}(1)}  \label{proofofsemisimplicity}

\sssec{}

To simplify notation we will assume that $\xl$ consists of one
point $x$. The proof in the general case is exactly the same.

To show that $\W^x_\varpi$ is semi-simple, it suffices to prove that
$$\on{Ext}^1(\Psib^{x,\mu}_\varpi,\Psib^{x,\nu}_\varpi)=0$$
for any $\mu,\nu \in \Lambda^{++}$.

\sssec{\bf Step 1}

We claim that it is enough to prove that
\begin{equation}    \label{ext0}
\on{Ext}^1(\Psib^{x,0}_\varpi,\Psib^{x,\nu}_\varpi)=0.
\end{equation}

Indeed, using \thmmref{two} and \propref{properties}(2,3) we obtain:

$$\on{Ext}^1(\Psib^{x,\mu}_\varpi,\Psib^{x,\la}_\varpi) \simeq
\on{Ext}^1(\Psib^{x,0}_\varpi\star\A_\mu,\Psib^{x,\la}_\varpi)
\simeq \on{Ext}^1(\Psib^{x,0}_\varpi, 
\Psib^{x,\la}_\varpi\star\A_{-w_0(\mu)}),$$ 
which is a direct sum of terms of the form
$\on{Ext}^1(\Psib^{x,0}_\varpi,\Psib^{x,\nu}_\varpi),
\nu \in \Lambda^{+}$, according to \corref{decomp}.

\sssec{\bf Step 2}

We claim that it suffices to show that for all $\la\in\Lambda^{+}$,
\begin{equation}    \label{last}
\on{Ext}^1(\Psib^{x,\la}_\varpi,\Psib^{x,\la}_\varpi)=0.
\end{equation}

First of all, \eqref{ext0} is obvious unless
$\nu\in\on{Span}\{\alpha_\i\}$, for otherwise $\Psib^{x,0}_\varpi$ and
$\Psib^{x,\nu}_\varpi$ live on different connected components of
$\bnx$.  Hence we can assume that $\nu\in\on{Span}\{\alpha_\i\}$,
which means that $Z(\GL)$ acts trivially on $V^\nu$.

The following result is well-known.

\bigskip

\noindent{\bf Lemma.}
{\em Let $H$ be a reductive algebraic group over an
algebraically closed field of characteristic $0$ and $V$ be an
irreducible $H$--module such that $Z(H)$ acts trivially on $V$. Then
there exists another irreducible $H$-module $W$, such that $V$ is a
direct summand in $W \otimes W^*$.}

\bigskip

This lemma, combined with \corref{decomp}, shows that for
$\nu\in\Lambda^{++}$, there exists $\la \in\Lambda^{++}$, such that
$\Psib^{x,\nu}_\varpi$ is a direct summand in
$\Psib^{x,\la}_\varpi\star \A_{-w_0(\la)}$.

Using Step 1, it is therefore enough to show that
$$\on{Ext}^1(\Psib^{x,0}_\varpi,\Psib^{x,\la}_\varpi\star \A_{-w_0(\la)})$$ 
for every $\la\in\Lambda^{++}$. But
\propref{properties}(2,3) and \thmmref{two} imply:
\begin{equation}
\on{Ext}^1(\Psib^{x,0}_\varpi,\Psib^{x,\la}_\varpi\star \A_{-w_0(\la)})\simeq
\on{Ext}^1(\Psib^{x,\la}_\varpi,\Psib^{x,\la}_\varpi).
\end{equation}

\sssec{\bf Step 3}

Finally, to prove \eqref{last} let us recall the following well-known
property of the Goresky-MacPherson extension.

\bigskip

\noindent{\bf Lemma.}
{\em Let $Y_0\overset{j}\hookrightarrow Y$ be an embedding
of algebraic stacks and let $\S_1$ and $\S_2$ be two perverse sheaves
on $Y_0$. Then the restriction map
$$\on{Ext}^1_Y(j_{!*}(\S_1),j_{!*}(\S_2))\to
\on{Ext}^1_{Y_0}(\S_1,\S_2)$$ is injective.}

\bigskip

Hence, it is enough to prove that
\begin{equation}    \label{ext1}
\on{Ext}^1_{_{\la,x}\Bun_N^{\F_T}}(\Psi^{x,\la}_\varpi,\Psi^{x,\la}_\varpi)=0.
\end{equation}

Recall from \secref{independence} that we can identify
$_{\la,x}\Bun_N^{\F_T}$ with $_{0,x}\Bun_N^{\F'_T} = \Bun_N^{\F'_T}$
for $\F'_T = \F_T(-\la \cdot x)$. Moreover, there exist $\varpi'_\i:
\L^{\alphach_\i}_{\F'_T}\to \Omega$ of conductor $\la \cdot x$, such
that under this identification the sheaf $\Psi^{x,\la}_\varpi$ becomes
$\Psi^{x,0}_{\varpi'} = \Psi_{\varpi'}$. Therefore \eqref{ext1} is
equivalent to
\begin{equation}    \label{ext2}
\on{Ext}^1_{\Bun_N^{\F'_T}}(\Psi_{\varpi'},\Psi_{\varpi'}) = 0.
\end{equation}

Since $\Psi_{\varpi'}$ is a rank $1$ local system on $\Bun_N^{\F'_T}$,
\eqref{ext2} is equivalent to
$$\on{Ext}^1_{\Bun_N^{\F'_T}}(\Ql,\Ql)=0.$$

According to  \corref{tower},  for $k$ large, $\Bun_N^{\F'_T} \simeq
{}^k\canN/N'_k$ and the space  $^k\canN_y$ is
isomorphic to a tower of affine spaces. Hence
$$
\on{Ext}^1_{\Bun_N^{\F'_T}}(\Ql,\Ql)\
=\ \on{Ext}^1_{^k\canN}(\Ql,\Ql)\ =\ H^1({}^k\canN)=0\,,
$$
where the first isomorphism follows from the fact that $N'_k$ is connected.

This completes the proof of \thmmref{one}(1).

\ssec{Proof of \thmmref{one}(2)}    \label{structure}

\sssec{}

The key step in proving part (2) of \thmmref{one} is the following
proposition.

\begin{prop}     \label{stability}
(1) The complex $\Psi^{\xl,\nul}_\varpi{}_{!}$ is a perverse sheaf.

(2) $\Psi^{\xl,\nul}_\varpi{}_{!}$ belongs to the category
$\W^{\xl}_\varpi$.
\end{prop}

Let $j: Y_0 \hookrightarrow Y$ be an affine embedding, and $\E$ a
local system on $Y_0$. Then the sheaf $j_! (\E)$ is perverse sheaf.

Therefore statement (1) of \propref{stability} follows from
\propref{affine}.

\medskip

Now we turn to part (2) of \propref{stability}.

\sssec{}
Note that $_{\xl,\nul}\ol{\Bun}_N^{\F_T}$ is the union of the strata
$_{\xl \cup \zl,\mul \cup \lal}\Bun_N^{\F_T}$, where $\mu_i\leq \nu_i,
i=1,\ldots,n$, $\zl = \{z_1,\ldots,z_m\}, \lal =
\{\la_1,\ldots,\la_m\}$, and $\la_j \in -\Lambda^+, j=1,\ldots,m$.

Call a stratum {\em relevant} if $\lal=0$ and $\mu_i\in\Lambda^{++}$
for all $i$ and irrelevant otherwise.  \propref{stability}(2) will be
derived from the following statement.

\sssec{\bf Lemma.}    \label{stalkvanishing}

\begin{enumerate}

\item[(1)] {\em The $*$--restriction of $\Psib^{\xl,\nul}_\varpi$ to
any irrelevant stratum is $0$.}

\item[(2)] {\em The perverse cohomologies of the $*$--restriction of
$\Psib^{\xl,\nul}_\varpi$ to a relevant stratum \newline
$_{\xl,\mul}\Bun_N^{\F_T}$ are direct sums of  copies of
 $\Psi^{\xl,\mul}_\varpi$.}

\end{enumerate}

\bigskip

Replacing $\F_T$ by
$$
\F'_T = \F_T \left( - \sum_{i=1}^n \nu_i \cdot x_i \right)
$$
transforms the sheaf $\Psib^{\xl,\nul}_\varpi$ over
$_{\xl,\nul}\ol{\Bun}_N^{\F_T}$ to $\Psib_{\varpi'}$ over
$\Bun_N^{\F'_T}$, where $\varpi'$ has conductor $\sum_i
\on{pr}(\nu_i)\cdot x_i$. Thus \lemref{stalkvanishing} can be
reformulated as follows:

\sssec{\bf Lemma.}  \label{sstalks} {\em The perverse cohomologies of the
$*$--restriction of
$\Psib_{\varpi'}$ to a stratum of the form $_{\xl,\mul}\Bun_N^{\F'_T}$
(where $\mu_i \in -\Lambda^+$) are zero unless
$\on{pr}(\mu_i)+\on{cond}_{x_i}(\varpi') \in
\Lambda_{G_{\on{ad}}}^{++}$ for every $i$. In the latter case, they are
direct sums of  copies of
$\Psi^{\xl,\mul}_{\varpi'}$.}

\sssec{} Let $y$ be a point different from the $x_i$'s and recall that
we have a sequence of inclusions:
$$_{\xl,\mul}\Bun_N^{\F'_T}\hookrightarrow
{}_{y,0}\wt{\Bun}_N^{\F'_T}\hookrightarrow \ol{\Bun}_N^{\F'_T}.$$

Recall also the $N(\O_y)$--bundles $\canN_y$ and $\CanN_y$ over
$\Bun_N^{\F'_T}$ and $_{y,0}\wt{\Bun}_N^{\F'_T}$, respectively. We will denote
the restriction of $\CanN_y$ to the stratum $_{\xl,\mul}\Bun_N^{\F'_T}$
by $_{\xl,\mul}\canN_y$. In addition, we will denote by $^k\canN_y$,
$^k\CanN_y$ and $_{\xl,\mul}^k\canN_y$ the corresponding
$N'_k$--bundles.

As shown in \secref{extension}, the stacks $^k\canN_y$, $^k\CanN_y$
and $_{\xl,\mul}^k\canN_y$ carry an action of a bigger group, namely,
$N_k$.  The proof of \lemref{sstalks} will be obtained by
analyzing this action.

\sssec{The additive character of $N(\K_y)$}  \label{groupcharacters}

The data of $\epsilon$ and $\varpi'$ give rise to a homomorphism
$\chi:N(\K_y) \to\GG_a$. Indeed, consider the ind--group scheme
$\Omega_{\K_y}=\Gamma(\D^\times_y,\Omega)$. The residue map gives rise
to a homomorphism $\on{Res}:\Omega_{\K_y}\to\GG_a$.
The additive character $\chi$ is defined as the sum
$$
\chi = \sum_{\i\in\I} {}^\i\chi,
$$
where ${}^\i\chi$ is defined as the composition
$$N(\K_y) \to N/[N,N](\K_y) \overset{\i}{\to}
\GG_a(\K_y)\overset{\varpi'\circ\epsilon}\longrightarrow
\Omega_{\K_y}\overset {\on{Res}}\longrightarrow \GG_a.$$

It is easy to see that $\chi$ vanishes on $N(\O_y)$. Therefore, for each
$k>0$, we obtain an additive character
$\chi^k: N_k\to\GG_a$.

\medskip

In what follows, given an $N_k$--stack $Y$ and a morphism $Y\to
\GG_a$, we will say that the latter is $(N_k,\chi^k)$--equivariant if
the diagram
$$
\CD
N_k\times Y  @>>> N_k\times \GG_a \\
@VVV @VVV \\
Y @>>> \GG_a
\endCD
$$
is commutative, where the right vertical arrow is the composition of
$\chi$ and the addition in $\GG_a$.

If $Y$ is an $N_k$--stack and $\S$ is an object of $\Sh(Y)$, we shall
say that $\S$ is $(N_k,\chi^k)$--equivariant if its pull-back under
the action map $N_k\times Y\to Y$ is isomorphic to \newline
$\chi^k{}^*(\J_\psi)\boxtimes \S$ so that the natural associativity
requirement holds (recall that $\J_\psi$ denotes the Artin-Schreier
sheaf on $\GG_a$).

Evidently, if $Y\to \GG_a$ is an $(N_k,\chi^k)$--equivariant morphism,
then the pull-back of $\J_\psi$ on $Y$ with respect to this map is
$(N_k,\chi^k)$--equivariant as a complex.

\sssec{The $N(\K_y)$--equivariance.}

The crucial observation is that the map
$$^k\canN_y\to \Bun_N^{\F'_T}\overset{\on{ev}_{\varpi'}}\longrightarrow
\GG_a$$ is $(N_k,\chi^k)$--equivariant. This follows from the definition of
$\chi^k$.

Hence, the pull-back of $\Psi_{\varpi'}$ under the natural projection
$^k\canN_y\to\Bun_N^{\F'_T}$ is an $(N_k,\chi^k)$--equivariant sheaf on
$^k\canN_y$. Therefore, by functoriality, the pull-back of
$\Psib_{\varpi'}$ to $^k\CanN_y$ and the perverse
cohomologies of the $*$--restriction of the latter to
$_{\xl,\mul}^k\canN_y$ are $(N_k,\chi^k)$--equivariant.

Therefore, \lemref{sstalks} is a consequence of the following general result
on equivariant sheaves on $_{\xl,\mul}^k\canN_y$:

\sssec{\bf Lemma.}  \label{ssstalks} {\em Let $\S$ be an
$(N_k,\chi^k)$--equivariant perverse sheaf on $_{\xl,\mul}^k\canN_y$.
Assume in addition that $k$ is large enough.  Then $\S$ vanishes unless
for every $i$, $\on{pr}(\mu_i) + \on{cond}_{x_i}(\varpi') \in
\Lambda_{G_{\on{ad}}}^{++}$. In the latter case, $\S$ is a direct sum
of copies of the  pull-back of
$\Psi_{\varpi'}^{\xl,\mul}$ from $_{\xl,\mul}\Bun_N^{\F'_T}$ to
$_{\xl,\mul}^k\canN_y$ (appropriately shifted).}

\sssec{Proof of \lemref{ssstalks}}

We know from \propref{transitivity} that for $k$ large enough
$_{\xl,\mul}^k\canN_y$ is a scheme and can be identified with the
quotient $N^{\F''_T}_{\out,k}\backslash N_k$, with
$\F''_T=\F'_T\left(-\underset{i=1}{\overset{n}\sum} \mu_i \cdot x_i
\right)$. In particular, $_{\xl,\mul}^k\canN_y$ is a homogenous space for
$N_k$.

Consider the case when $\on{pr}(\mu_i) + \on{cond}_{x_i}(\varpi') \in
\Lambda_{G_{\on{ad}}}^{++}$. The
pull-back of
$\Psi^{\xl,\mul}_{\varpi'}$ to $_{\xl,\mul}^k\canN_y$
is $(N_k,\chi^k)$--equivariant; we see that using the same reasoning as
for $\Psi_{\varpi'}$. Since the isotropy subgroup
$N^{\F''_T}_{\out,k}$ is connected, any $(N_k,\chi^k)$--equivariant
perverse sheaf is isomorphic to a direct sum of copies of this pull-back
(up to the appropriate shift).

To prove the vanishing part of the statement it suffices to show the following: if
$\on{pr}(\mu_i)+\on{cond}_{x_i}(\varpi')$ is non-dominant for some $i$
then the homomorphism
$$N^{\F''_T}_{\out,k}\hookrightarrow N_k\overset{\chi^k}\longrightarrow
\GG_a$$ is non-trivial, for $k$ sufficiently large. Clearly, it
suffices to show that the composition $N^{\F''_T}_{\out}\to
N(\K_y)\overset{\chi}\to \GG_a$ is non-zero.

By definition, the above homomorphism $N^{\F''_T}_{\out}\to\GG_a$ is the
composition:
$$N^{\F''_T}_{\out} \twoheadrightarrow \prod_{\i\in\I}
H^0(X\setminus y,\L_{\F''_T}^{\alphach_\i}) \to
\prod_{\i\in\I} H^0(\D_y^\times, \L_{\F'_T}^{\alphach_\i})
\overset{\varpi'}\longrightarrow
\Omega_{\K_y}^\I
\overset{\on{Res}}\longrightarrow \GG_a^\I
\overset{\on{sum}}\longrightarrow \GG_a.$$

Assume that for some $j\in\{1,...,n\}$ and $\i\in\I$,
$\langle\on{pr}(\mu_j)+\on{cond}_{x_j}(\varpi'),\check\alpha_\i\rangle<
0$.  Then there exists a section $\gamma\in H^0(X\setminus
y,\L_{\F''_T}^{\alphach_\i})$ with the following properties:

\begin{itemize}

\item
$\gamma$ does not vanish under the composition
$$H^0(X\setminus y,\L_{\F''_T}^{\alphach_\i})\to
H^0(\D^\times_{x_j},\L_{\F''_T}^{\alphach_\i}) \overset{\varpi'_\i}\to
\Omega_{\K_{x_j}} \overset{\on{Res}}\longrightarrow\GG_a.$$

\item
$\gamma$ vanishes to a sufficiently high order at the points $x_{j'},
j'\neq j$, so that it goes to $0$ under the composition
$$H^0(X\setminus y,\L_{\F''_T}^{\alphach_\i})\to
H^0(\D^\times_{x_{j'}},\L_{\F''_T}^{\alphach_\i}) \overset{\varpi'_\i}\to
\Omega_{\K_{x_{j'}}} \overset{\on{Res}}\longrightarrow\GG_a.$$

\end{itemize}

But then by the residue formula, the image of $\gamma$ under the
composition
$$H^0(X\setminus y,\L_{\F''_T}^{\alphach_\i}) \to
H^0(\D_y^\times, \L_{\F'_T}^{\alphach_\i})
\overset{\varpi'_\i}\longrightarrow
\Omega_{\K_y}
\overset{\on{Res}}\longrightarrow \GG_a$$
is non-zero. Therefore $\chi(\gamma) \neq 0$. This shows that $\S=0$
under the above assumption.

This completes the proof of \lemref{ssstalks} and hence of
\lemref{stalkvanishing} and \lemref{sstalks}.

\sssec{}

In order to complete the proof of \propref{stability}(2), we need to
show that the sheaf $\Psi^{\xl,\nul}_\varpi{}_{!}$ has a filtration
whose consecutive quotients are perverse sheaves of the form
$\Psib^{\xl,\mul}_\varpi$.  For simplicity, consider the case $n=1$.
We prove this statement by induction on
$\langle\nu,\rhoc\rangle$.

If $\nu = 0$, then there are no relevant strata in the closure of
$_{x,\nu}\Bun_N^{\F_T}$ and therefore $\Psib^{x,0}_\varpi =
\Psi^{x,0}_\varpi{}_{!}$.

Now suppose that we have already proved the result for all dominant
weights $\nu$, such that $\langle\nu,\rhoc\rangle < N$. Consider a
dominant weight $\nu$, such that $\langle\nu,\rhoc\rangle=N$.

According to \propref{stability}(1), $\Psi^{x,\nu}_\varpi{}_{!}$ is a
perverse sheaf. Let ${\mathcal K}$ denote the perverse sheaf that 
fits into the exact sequence
$$
0 \to {\mathcal K} \to \Psi^{x,\nu}_\varpi{}_{!} 
\to \Psib^{x,\nu}_\varpi \to 0.
$$
But then ${\mathcal K}$ is non-zero only on the relevant strata
$_{x,\mu}\Bun_N^{\F_T}$ with $\langle\mu,\rhoc\rangle < N$, and we
know that the perverse cohomology groups of its restrictions to those
strata are direct sums of copies of
$\Psi^{x,\mu}_\varpi$. Hence, on the level of Grothendieck groups, 
we can write ${\mathcal K}$
as a combination of the sheaves $\Psi^{x,\mu}_\varpi{}_{!}$ with
$\langle\mu,\rhoc\rangle < N$. By our inductive assumption, the sheaves
$\Psi^{x,\mu}_\varpi{}_{!}$ are consecutive extensions of  perverse
sheaves of the form $\Psib^{x,\mu'}_\varpi$, and hence so is 
${\mathcal K}$ and, finally, so is $\Psi^{x,\nu}_\varpi{}_{!}$.

Thus, \propref{stability} is proved.

\sssec{Proof of \thmmref{one}(2)}

Let us show that the map
\begin{equation}    \label{mapp}
\Psi^{\xl,\nul}_\varpi{}_{!}\twoheadrightarrow
\Psib^{\xl,\nul}_{\varpi}
\end{equation}
is an isomorphism. The isomorphism
$\Psi^{\xl,\nul}_\varpi{}_{*} \simeq
\Psib^{\xl,\nul}_{\varpi}$
will then also follow by Verdier duality.

Let ${\mathcal K}$ denote the kernel of the map \eqref{mapp}. According to
\propref{stability}, ${\mathcal K}$ is an object of $\W^{\xl}_\varpi$ and by
\thmmref{one}(1), we have:
$$\Psi^{\xl,\nul}_\varpi{}_{!}\simeq 
\Psib^{\xl,\nul}_{\varpi}\oplus {\mathcal K}.$$

However, the restriction of $\Psi^{\xl,\nul}_\varpi{}_{!}$ to the
complement of $_{\xl,\nul}\BunNft$ in $_{\xl,\nul}\BunNftb$ is zero,
by definition. Therefore the same is true for ${\mathcal K}$.

But then ${\mathcal K}=0$, since by construction ${\mathcal K}$ is supported 
on $_{\xl,\nul}\BunNftb\setminus {}_{\xl,\nul}\BunNft$. Hence
\eqref{mapp} is an isomorphism.

\section{Computation of the Hecke action} \label{heckecontd}

In this section we collect several facts concerning the geometry of
$N(\K_x)$--orbits on $\Gr$ and prove \thmmref{two}.

\ssec{Semi-infinite orbits on $\Gr$.}

\sssec{}  \label{introrbits}

Recall the ind--scheme $\Gr$. Informally, one defines the locally
closed ind--subscheme $S^\nu\subset\Gr$ as the $N(\K_x)$--orbit of the
point $\nu(t)\in\Gr$. Here is a precise scheme--theoretic definition.

The chosen Borel subgroup $B\subset G$ defines for each
$\check\la\in\check\Lambda^{++}$ a line subbundle $\OO\subset
\V^{\check\la}_{\F^0_G}$ over $\D_x$. We say that a point
$(\F_G,\beta)\in\Gr$ belongs to $\ol{S}^\nu$, if for every $\check\la$
the map $\OO\to \V^{\check\la}_{\F_G}(\infty\cdot x)$ obtained using
the data of $\beta$ factors as
\begin{equation}    \label{factorsas}
\OO\to \V^{\check\la}_{\F_G}(\langle\nu,\check\la\rangle\cdot x)
\hookrightarrow \V^{\check\la}_{\F_G}(\infty\cdot x).
\end{equation}

It follows from the definition that there is a morphism $\ol{S}^\nu
\to {}_{x,\nu}\ol{\Bun}_N^{\F_T^0}$, which corresponds to forgetting
$\beta$, but keeping the maps $\OO \hookrightarrow
\V^{\check{\lambda}}_{\F_G}(\langle \nu,\check{\lambda} \rangle \cdot
x)$ induced by $\beta$. One can show that this morphism is formally
smooth.

It is clear that $\ol{S}^\nu$ is a closed ind-subscheme of $\Gr$ and
$\ol{S}^{\nu'}\subset \ol{S}^\nu$ if and only if $\nu\geq\nu'$. We
define the locally closed ind-subscheme $S^\nu\subset\Gr$ as follows:
we say that $(\F_G,\beta)\in\Gr$ belongs to $S^\nu$, if for every
$\check\la$ the map $\OO\to \V^{\check\la}_{\F_G}(\infty\cdot x)$
factors as \eqref{factorsas} and the embedding $\OO\to
\V^{\check\la}_{\F_G}(\langle\nu,\check\la\rangle\cdot x)$ is maximal.

The group $N(\K_x)$ acts naturally on $\Gr$. It is easy to see that each
subscheme $S^\nu$ is $N(\K_x)$--stable. Moreover, one can show
that the $N(\K_x)$--action on $S^\nu$ is transitive in the sense that
$S^\nu=\underset{k}\cup (S^\nu)^k$, where each $(S^\nu)^k$ is a
homogeneous space for $N^k$.

\sssec{Remark.}    \label{nilpotents}
It is instructive to compute explicitly $R$-points of the scheme
$\ol{S}^\nu$, when $G=SL_2$ and $R$ is an Artinian local ring. 
In this case, the set $\ol{S}^n(R) \subset \Gr(R)$ consists of the 
cosets of the form
$$
\begin{pmatrix} 1 & u(t) \\ 0 & 1 \end{pmatrix} \begin{pmatrix}
t^{n-m} p_m(t^{-1})^{-1} & 0 \\ 0 & t^{-n+m} p_m(t^{-1}) \end{pmatrix}
\cdot SL_2(R[[t]]),\,\, m\geq 0
$$
where $u(t) \in R((t))$, and $p_m(t^{-1})$ is an invertible element of
$R[t^{-1}]$ of degree $m$. This means that (up to a scalar in
$R^\times$, which can be absorbed into $SL_2(R[[t]])$)
$$
p_m(t^{-1}) = 1 + r_1 t^{-1} + \ldots + r_m t^{-m},
$$
where each coefficient $r_i$ is a nilpotent element of $R$.

Thus, we see that the set of $R$--points of $\ol{S}^n$ coincides with
the set of $R$--points of a union of strata which are fibrations over
$Sym^m \D_x$ with fibers $S^{n-m}$, where $m \geq 0$.
\footnote{Here $Sym^m \D_x$ is a formal scheme
$\on{Spf}(\Fqb[[t_1,...,t_m]]^{\Sigma^m})$, where $\Sigma^m$ is the
symmetric group.} These fibrations are similar to those described in
\corref{stratified}: the formal disc $\D_x$ here plays the role of the
curve $X$.

This is analogous to the description of $\ol{S}^n$ obtained in
\cite{FS} in the case when $\O_x$ is replaced by the ring of analytic
functions on the unit disc -- then $\D_x$ is replaced by the unit disc.

\sssec{} \label{dimestimate}

The following proposition, due to \cite{MV}, will play an important
role in the proof of \thmmref{two}.

\medskip

\noindent{\bf Proposition.}
$\on{dim}(\Gr^\la\cap S^\nu)\leq \langle \la+\nu,\rhoc\rangle$.

\bigskip

It is also known that for $\nu=\la$, $\Gr^\la\cap S^\nu$ is open and
dense in $\Gr^\la$ and for $\nu=w_0(\la)$, $\Gr^\la\cap S^\nu$ is a
point--scheme.

\sssec{Admissible characters}  \label{admcharacters}

Let $\eta$ be a coweight of $G_{\on{ad}}$ and let us choose
isomorphisms of $\O_x$--modules
\begin{equation}  \label{datachar}
\O_x(\langle\eta,\check\alpha_\i\rangle \cdot x)\simeq \Omega|_{\D_x}
\end{equation}
for each $\i\in\I$.

As in \secref{groupcharacters}, the above data define a homomorphism
$\chi_\eta:N(\K_x)\to\GG_a$. We call such a character
admissible of conductor $\eta$. It is clear that the
group $T(\O_x)(\Fqb)$ acts transitively on the set admissible
characters with a given conductor.

In order to simplify notation, we will write $\chi_\nu$ for
$\chi_{\on{pr}(\nu)}$ for any $\nu \in \Lambda$.

\sssec{\bf Lemma.}  \label{orbitchar} {\em Let $\nu$ and $\mu$ be two
elements of $\Lambda$. Then for a given admissible character
$\chi_\mu$ of conductor $\on{pr}(\mu)$ there exists a
$(N(\K_x),\chi_{\mu})$--equivariant function $\chi^\nu_{\mu}: S^\nu
\to \GG_a$ if and only if $\mu+\nu \in \Lambda^{++}$.  In the latter
case this function is unique up to an additive constant.}

\sssec{}

When $\mu+\nu \in \Lambda^{++}$, the  function $\chi^\nu_\mu: S^\nu \to
\GG_a$ can be defined by the formula
\begin{equation}    \label{chimunu}
\chi^\nu_\mu(n \cdot \nu(t)) = \chi(\mu(t)\cdot n \cdot\mu^{-1}(t)), \quad
\quad n \in N(\K).
\end{equation}

The next statement, which is proved in \secref{proofcomponents}, is
a key ingredient in the proof of \thmmref{two}.

\sssec{\bf Proposition.} \label{components} {\em Let
$\chi^{\nu}_{-\nu}: S^{\nu}\to \GG_a$ be an
$(N(\K_x),\chi_{-\nu})$--equivariant function. Assume that
$\la\in\Lambda$ is a dominant coweight and $\nu \neq w_0(\la)$. Then
the restriction of $\chi^{\nu}_{-\nu}$ to any irreducible
component of $\Gr^\la\cap S^{\nu}$ of dimension $\langle
\la+\nu,\rhoc\rangle$ is a dominant map onto $\GG_a$.}

\sssec{\bf Corollary.} \label{corcomp}
{\em The cohomology
$H^{\langle\la+\nu,2\rhoc\rangle}_c(\Gr^\la\cap S^{\nu},
\chi^{\nu}_{-\nu}|_{\Gr^\la\cap S^{\nu}}{}^*(\J_\psi))$
vanishes unless $\nu=w_0(\la)$.}

\sssec{Proof.}

Without loss of generality we can assume that
$\chi^{\nu}_{-\nu}(\nu(t))=0\in\GG_a$.

Denote by $K$ the complex on $\GG_a$, which is the !--direct image of the
constant sheaf on $\Gr^\la\cap S^\nu$ with respect to the map
$\chi^{\nu}_{-\nu}: \Gr^\la\cap S^{\nu}\to \GG_a$.

Consider the $T$--action on $\GG_a$ given by the character
$2\rhoc:T\to\GG_m$. It is easy to see that $K$ is $T$--equivariant
with respect to this action. In particular, it is monodromic with
respect to the standard $\GG_m$--action on $\GG_a$.

\propref{components} and \propref{dimestimate} imply that our complex
lives in the (perverse) cohomological degrees strictly less than
$\langle\la+\nu,2\rhoc\rangle$. \corref{corcomp} is therefore
equivalent to saying that $H_c^1(\GG_a,\wt{K}\otimes \J_\psi)=0$, where
$\wt{K}$ is the $(\langle\la+\nu,2\rhoc\rangle-1)$-th perverse
cohomology sheaf of $K$. But this follows from:

\sssec{\bf Lemma.}  \label{AS} {\em Let $\S$ be a $\GG_m$--monodromic
perverse sheaf on $\AA^1$. Then we have: $H_c^1(\AA^1,\S\otimes
\J_\psi)=0$.}

\ssec{Proof of \thmmref{two}.}

\sssec{} For notational convenience, we replace $\lambda$ by $-w_0(\la)$
and prove the statement:
$$\Psib_\varpi^{x,0} \star \A_{-w_0(\la)}=
\Psib_\varpi^{x,-w_0(\la)}\,.$$

Thus, we need to establish an isomorphism
\begin{equation} \label{repeatclaim}
'\hl_{!}((\Psib_\varpi^{x,0}\tboxtimes\A_{-w_0(\la)})^r) \simeq
\Psib_\varpi^{x,-w_0(\la)}.
\end{equation}
To simplify notation, from now on we will suppress the upper index $r$
in this formula.

Denote by $\wt{K}^\nu$ (resp., $K^\nu$) the $*$--restriction
of the LHS of \eqref{repeatclaim} to the stratum $_{x,\nu}\BunNftw$
(resp., $_{x,\nu}\BunNft$). Since the both sides  of \eqref{repeatclaim}
are Verdier self--dual (up to replacing $\psi$ by $\psi^{-1}$) it suffices
to prove  the following

\begin{claimm} \label{strat1}
\begin{itemize}

\item[(1)] The complex $\wt{K}^\nu$ lives in the (perverse)
cohomological degrees $\leq 0$.

\item[(2)]
The $*$--restriction of $\wt{K}^{\nu}$ to the closed substack
$$_{x,\nu}\BunNftw\setminus {}_{x,\nu}\BunNft\subset {}_{x,\nu}\BunNftw$$
lives in strictly negative (perverse) cohomological degrees.

\item[(3)]
The $0$-th (perverse) cohomology of $K^{\nu}$ vanishes
unless $\nu=-w_0(\la)$ and in the latter case it can be identified
with $\Psi^{x,-w_0(\lambda)}_\varpi$.

\end{itemize}
\end{claimm}

\sssec{}  \label{stratifications}

For each pair $\nu,\nu'\in\Lambda$ consider the following locally closed
substacks of $Z$:
\begin{align*}
&\wt{Z}^{\nu,?}:={}'\hl{}^{-1}({}_{x,\nu}\BunNftw),\,\,\,
Z^{\nu,?}:={}'\hl{}^{-1}({}_{x,\nu}\BunNft) \\
&\wt{Z}^{?,\nu'}:={}'\hr{}^{-1}({}_{x,\nu'}\BunNftw),\,\,\,
Z^{?,\nu'}:={}'\hr{}^{-1}({}_{x,\nu'}\BunNft) \\
&\wt{Z}^{\nu,\nu'}:=\wt{Z}^{\nu,?}\cap \wt{Z}^{?,\nu'},\,\,\,
Z^{\nu,\nu'}:=Z^{\nu,?}\cap \wt{Z}^{?,\nu'}.
\end{align*}

For $\mu\in\Lambda^{++}$, denote by $\H_x^\mu$ the locally closed
substack $\canG_x\underset{G(\O_x)}\times\Gr^\mu$ of $\H_x$
(the projection onto $\Bun_G$ that we are using here is $\hl$), and set

\begin{align*}
&\wt{Z}^{\nu,?,\mu}:=\wt{Z}^{\nu,?}\cap {}'\p^{-1}(\H_x^\mu),\,\,\,
\wt{Z}^{?,\nu',\mu}:=\wt{Z}^{?,\nu'}\cap {}'\p^{-1}(\H_x^\mu) \\
&\wt{Z}^{\nu,\nu',\mu}:=\wt{Z}^{\nu,\nu'}\cap {}'\p^{-1}(\H_x^\mu),\,\,\,
Z^{\nu,\nu',\mu}:=Z^{\nu,\nu'}\cap
{}'\p^{-1}(\H_x^\mu)
\end{align*}

Let us denote by $\wt{K}^{\nu,\nu',\mu}$ (resp., by
$K^{\nu,\nu',\mu}$) the $!$-direct image under
$$'\hl: \wt{Z}^{\nu,\nu',\mu}\to {}_{x,\nu}\BunNftw \text{ (resp., }
Z^{\nu,\nu',\mu}\to {}_{x,\nu}\BunNft \text{ ) }$$
of the $*$--restriction of $\Psib_\varpi^{x,0}\tboxtimes\A_{-w_0(\la)}$
to $\wt{Z}^{\nu,\nu',\mu}$ (resp., to $Z^{\nu,\nu',\mu}$).

Using a standard spectral sequence, one can derive
\claimmref{strat1} from the following

\begin{claimm} \label{strat2}

\begin{itemize}
\item[(1)]
The complex $\wt{K}^{\nu,\nu',\mu}$ lives in  cohomological degrees
$\leq 0$ and the inequality is strict unless $\nu'=0$ and $\mu=\la$.

\item[(2)]
The $*$--restriction of $\wt{K}^{\nu,0,\la}$ to the closed substack
$$_{x,\nu}\BunNftw\setminus {}_{x,\nu}\BunNft\subset {}_{x,\nu}\BunNftw$$
lives in strictly negative cohomological degrees.

\item[(3)]
The $0$-th cohomology of $K^{\nu,0,\la}$
vanishes unless $\nu=-w_0(\la)$.

\item[(4)]
The $0$-th cohomology of $K^{-w_0(\la),0,\la}$
can be identified with $\Psi^{x,-w_0(\lambda)}_\varpi$.

\end{itemize}

\end{claimm}

The proof of \claimmref{strat2} will be obtained by analyzing the
fibers of the projection
$'\hl: \wt{Z}^{\nu,\nu',\mu}\to {}_{x,\nu}\BunNftw$.

\sssec{}

Recall that by choosing a trivialization $\epsilon:\F_T\to
\F^0_T|_{\D_x}$, we obtain canonical $N(\O_x)$--torsors $\CanN_x$ and
$\canN_x$ over $_{x,0}\widetilde{\Bun}^{\F_T}_N$ and $\BunNft$,
respectively.

Independently, for every $\nu\in\Lambda$ and $\F'_T:=\F_T(-\nu\cdot
x)$ we fix a trivialization
$\epsilon_\nu:\F'_T\to\F^0_T|_{\D_x}$. Using the isomorphism
$_{x,\nu}\BunNftb\simeq \ol{\Bun}_N^{\F'_T}$ we obtain
$N(\O_x)$--torsors $_\nu\CanNnu_x$ and $_\nu\canNnu_x$ over
$_{x,\nu}\widetilde{\Bun}^{\F_T}_N$ and $_{x,\nu}\BunNft$,
respectively.

By definition, the projection $'\hl$ identifies $\wt{Z}^{\nu,?}$ with the
fibration
$$_\nu\CanNnu_x\underset{N(\O_x)}\times\Gr \to {}_{x,\nu}\BunNftw,$$
and similarly, the projection $'\hr$ realizes $\wt{Z}^{?,\nu'}$ as a
fibration
$$_{\nu'}\CanNnup_x\underset{N(\O_x)}\times\Gr \to {}_{x,\nu'}\BunNftw.$$

\medskip

Recalling the definition of the subscheme $S^\nu\subset\Gr$ from
\secref{introrbits}, we obtain the following lemma:

\sssec{\bf Lemma.}  \label{partdescr}

\begin{enumerate}

\item[(1)] {\em The stacks $\wt{Z}^{\nu,\nu'}$ and
$\wt{Z}^{\nu,?,\mu}$, when viewed as substacks of $\wt{Z}^{\nu,?}$,
can be identified with the fibrations}
$$_\nu\CanNnu_x\underset{N(\O_x)}\times
S^{\nu'-\nu}\overset{'\hl}\longrightarrow{}_{x,\nu}\BunNftw \text{ and
} _\nu\CanNnu_x\underset{N(\O_x)}\times
\Gr^{\mu}\overset{'\hl}\longrightarrow{}_{x,\nu}\BunNftw, \text{
respectively.}$$

\item[(2)] {\em The stacks $\wt{Z}^{\nu,\nu'}$ and
$\wt{Z}^{?,\nu',\mu}$, when viewed as substacks of $\wt{Z}^{?,\nu'}$,
can be identified with the fibrations}
$$_{\nu'}\CanNnup_x\underset{N(\O_x)}\times
S^{\nu-\nu'}\overset{'\hr}\longrightarrow{}_{x,\nu'}\BunNftw \text{
and } _{\nu'}\CanNnup_x\underset{N(\O_x)}\times
\Gr^{-w_0(\mu)}\overset{'\hr}\longrightarrow{}_{x,\nu'}\BunNftw,
\text{ respectively.}$$

\item[(3)] {\em The substack $Z^{\nu,\nu'}$ coincides with the
intersection $\wt{Z}^{\nu,?}\cap Z^{?,\nu'}$.}

\end{enumerate}

\sssec{}
Now, we are ready to prove the first two statements of \claimmref{strat2}.

Indeed, from \lemref{partdescr} we conclude that the $*$--restriction
of $\Psib_\varpi^{x,0}\tboxtimes\A_{-w_0(\la)}$ to $\wt{Z}^{?,\nu',\mu}$ is
the twisted external product of complexes:
\begin{equation}  \label{descrrestr}
(\Psib_\varpi^{x,0}|_{_{x,\nu'}\BunNftw})\tboxtimes
(\A_{-w_0(\la)}|_{\Gr^{-w_0(\mu)}}).
\end{equation}

By definition, it lives in the cohomological degrees $\leq 0$ and the
inequality is strict unless $\mu=-w_0(\la)$ and $\nu'=0$. (In
fact, we know from \propref{stability} that the above complex is $0$
unless $\nu'=0$, but we will not need this fact here.)

Now since $\A_{-w_0(\la)}|_{\Gr^{-w_0(\mu)}}$ is  constant,
it follows from \lemref{partdescr}(2) that its further $*$--restriction
to $\wt{Z}^{\nu,\nu',\mu}$ lives in (perverse) cohomological degrees
$$\leq
-\on{codim}(S^{\nu-\nu'}\cap\Gr^{-w_0(\mu)},\Gr^{-w_0(\mu)})\leq
-\langle -w_0(\mu)-\nu+\nu',\rhoc\rangle=
-\langle\mu-\nu+\nu',\rhoc\rangle,$$ according to
\propref{dimestimate}.

However, we obtain from \lemref{partdescr}(1) and
\propref{dimestimate} that the fibers of the projection $'\hl:
\wt{Z}^{\nu,\nu',\mu}\to{}_{x,\nu}\BunNftw$ are of dimension $\leq
\langle\mu-\nu+\nu',\rhoc\rangle$.

This proves \claimmref{strat2}(1).

\claimmref{strat2}(2) follows from \lemref{partdescr}(3) combined with
the above dimension estimates and the fact that the $*$--restriction
of $\Psib_\varpi^{x,0}$ to $_{x,0}\BunNftw\setminus {}_{x,0}\BunNft$
lives in strictly negative cohomological degrees.

Thus, it remains to study $K^{\nu,0,\la}$.

\sssec{}

Let $\canNK$ denote the ind--group scheme over
$_{x,\nu}\BunNft$, obtained as a $_\nu\canNnu_x$--twist of $N(\K_x)$
with respect to the adjoint action of $N(\O_x)$ on $N(\K_x)$. The
description of $Z^{\nu,\nu'}$ as a stack fibered over
$_{x,\nu}\BunNft$ in \lemref{partdescr}(1) implies that
$Z^{\nu,\nu'}$ carries a canonical
$\canNK$--action.

Recall now from \secref{admcharacters} that the data of
$\epsilon_\nu$ and $\varpi$ give rise to an admissible character
$\chi_{\nu}:N(\K_x)\to\GG_a$, of conductor $\on{pr}(\nu)$. Hence we 
also obtain a character on the above ind--group scheme
$\canNK$, which we denote by $\ol{\chi}_{\nu}$.

For dominant coweights $\nu$ and $\nu'$, the
$(N(\K_x),\chi_{\nu})$--equivariant function
$\chi_{\nu}^{\nu'-\nu}:S^{\nu'-\nu}\to\GG_a$ gives rise to a
$(\canNK,\ol{\chi}_{\nu})$--equivariant function
$\ol{\chi}_{\nu}^{\nu'-\nu}:Z^{\nu,\nu'}\to\GG_a$.

The following result is obtained directly from the definitions:

\sssec{\bf Lemma.}  \label{functionsmatch}
{\em Assume that $\nu'\in\Lambda$ is dominant. Then}

\begin{enumerate}

\item[(1)] {\em The function
$$(\on{ev}^{x,\nu'}_{\varpi}\circ{}'\hr):
Z^{\nu,\nu'}\to\GG_a$$ is
$(\canNK,\ol{\chi}_{\nu})$--equivariant.}

\item[(2)] {\em If $\nu$ is also dominant, then the above function
coincides with the composition
$$Z^{\nu,\nu'}\overset{\ol{\chi}^{\nu'-\nu}_{\nu}
\times{}'\hl}\longrightarrow \GG_a\times {}_{x,\nu}\BunNft
\overset{\on{id}\times\on{ev}^{\nu,x}_\varpi}\longrightarrow
\GG_a\times\GG_a\overset{\on{sum}} \longrightarrow \GG_a$$ for some
$\chi^{\nu'-\nu}_{\nu}$.}

\end{enumerate}

\sssec{} Asume that $\nu$ is different from $-w_0(\la)$.
\lemref{partdescr}(1) implies that after a
smooth localization $U\to {}_{x,\nu}\BunNft$ (e.g., we can take
$U={}_\nu^k\canNnu_x$, for large enough $k$), the fibration $'\hl:
Z^{\nu,0,\la}\to{}_{x,\nu}\BunNft$ becomes a direct product $U\times
(\Gr^\la\cap S^{-\nu})$.

Moreover, by \lemref{functionsmatch}(1), the complex
$(\Psib_\varpi^{x,0}\tboxtimes\A_{-w_0(\la)})|_{Z^{\nu,0,\la}}$
becomes the external product
$$\E\boxtimes \chi^{-\nu}_{\nu}{}^*(\J_\psi)[\langle\la-\nu,2\rhoc\rangle],$$
where $\E$ is a locally constant perverse sheaf on the base $U$.
Therefore \claimmref{strat2}(3) follows from \corref{corcomp}.

Finally, for $\nu=-w_0(\la)$ the intersection
$S^{w_0(\la)}\cap\Gr^\la$ is a point-scheme (cf. \secref{introrbits})
and therefore \claimmref{strat2}(4) follows from
\lemref{functionsmatch}(2) applied to the case $\nu'=0$.

In order to finish the proof of \claimmref{strat2} (and hence of
\thmmref{two}), it remains to prove \propref{components}.

\ssec{Proof of \propref{components}} \label{proofcomponents}

In order to prove the proposition, it suffices to study the function
$\chi^\nu_{-\nu}$ on the set of $\Fqb$--points of $\Gr^\la \cap S^\nu$.

\sssec{Computation in rank $1$}  \label{rk1}

First we will consider the case when $G$ has semi--simple rank $1$. We
can identify each intersection of an $N(\K_x)$--orbit and an
$G(\O_x)$--orbit in $\Gr_G$ with an appropriate intersection of an
$N(\K_x)$--orbit and an $G_{\on{ad}}(\O_x)$--orbit in $\Gr_{G_{\on{ad}}}$,
where $G_{\on{ad}} = G/Z(G)$. Moreover, the corresponding functions
$\chi_{-\nu}^{\nu}$ coincide under this identification. Therefore
without loss of generality we can replace $G$ by the corresponding
adjoint group, so it suffices to treat the case of the group $PGL(2)$.
In this case we prove the statement of \propref{components} by
an explicit computation as follows.

Let us identify $\Lambda$ with $\Z$.
The intersection $\ol{\Gr}^m \cap S^n$ is empty unless $m$ and $n$ have
the same parity and $|n| \leq m$.  When this is the case, $\ol{\Gr}^m
\cap S^n$ is isomorphic to ${\mathbb A}^{(n+m)/2}$:
\begin{equation}    \label{inters}
\ol{\Gr}^m \cap S^n = \left\{ \begin{pmatrix} 1 &
\sum_{i=(n-m)/2}^{n-1} a_i t^i \\ 0 & 1
\end{pmatrix} \cdot
\begin{pmatrix} t^n & 0 \\ 0 & 1 \end{pmatrix} \mid
a_i\in\Fqb\right\}\,.
\end{equation}
In particular, $\ol{\Gr}^m \cap S^{n}$ is
always irreducible.

Hence in order to prove \propref{components} it is enough to show that
on the level of $\Fqb$--points, the function
$$\chi^{n}_{-n}:(\ol{\Gr}^m \cap S^{n}) \to\Fqb$$ is non--constant
whenever $\ol{\Gr}^m \cap S^{n}$ has positive dimension.

Using the $a_i$ in \eqref{inters} as coordinates on $\ol{\Gr}^m
\cap S^n$, the function
$\chi_{-n}^n$ is given by $\psi(a_{n-1})$. Hence the restriction of
$\chi^{n}_{-n}$ to $\ol{\Gr}^m \cap S^n$ is indeed non-constant if the
latter has positive dimension.

\sssec{Grassmannians associated to the parabolic subgroups of $G$.}
The proof for general $G$ will be obtained by reduction from
$G$ to its minimal Levi subgroups.  In what follows, the subscript
``$M$'' (e.g. in $\Gr_M$, $S^\nu_M$, etc.) will denote the
corresponding object for a Levi subgroup $M$ of $G$.

Let $P$ be a parabolic subgroup of $G$, $N_P$ its unipotent radical, and
$M$ its Levi subgroup.  Denote by $\I_M \subset \I$ the 
coresponding subset of $\I$. Let
$\Lambda_{G,P}$ be the quotient of $\Lambda=\Lambda_M$ by the
sublattice spanned by $\al_{\i}, \i \in \I_M$. Then $\Lambda_{G,P}$
can be identified with the set of connected components of the
Grassmannian $\Gr_M = M(\K_x)/M(\O_x)$. We denote by $\Gr_M^{(\theta)}$ the
component of $\Gr_M$ corresponding to $\theta \in \Lambda_{G,P}$. Note
that $\Gr_M^\mu\subset\Gr_M^{(\theta)}$ if and only if
$\mu$ belongs to the coset of $\theta$ in $\Lambda$.

Let $P_0(\K_x)$ be a ind--subgroup of $P(\K_x)$, defined as the inverse
image of the subgroup
$M/[M,M](\O_x)$ under the natural projection $P(\K_x)\to M(\K_x)\to
M/[M,M](\K_x)$. As in \secref{introrbits}, to an element
$\theta\in\Lambda_{G,P}$ one can attach a locally closed
ind--subscheme $S^\theta_P$ in $\Gr$, on which $P_0(\K_x)$ acts
transitively (at the level of points, $S^\theta_P$ is the
$P_0(\K_x)$--orbit of $\theta(t)$ in $\Gr$).

In the same way as in \secref{affgrass} we can define the affine
Grassmannian $\Gr_P$ (note that the definition given in
\secref{affgrass} works for any algebraic group). The set of connected
components of $\Gr_P$ coincides with that of $\Gr_M$. Now
$S^\theta_P$ is nothing but the reduced scheme of the corresponding
connected component of $\Gr_P$. Thus, the natural projection
$\Gr_P\to\Gr_M$ gives rise to a map
$p^{(\theta)}: S_P^\theta\to\Gr_M^{(\theta)}$.

The following statement is straightforward.

\sssec{\bf Lemma.}
{\em For $\nu\in\Lambda$ let $\theta$ be its image in
$\Lambda_{G,P}$. Then $S^\nu\subset S^\theta_P$ and
$p^{(\theta)}(S^\nu)\subset S_M^\nu$.}

\bigskip

Denote by $p^\nu$ the restriction of $p^{(\theta)}$ to $S^\nu$.

\sssec{}

For $\i \in \I$, denote by $P_{\i}$ the corresponding minimal
parabolic subgroup of $G$ and by $M_{\i}$ its Levi subgroup.  Then we
have a morphism $p_{\i}^\nu: S^\nu \arr S_{M_{\i}}^\nu$. Note that
$S_{M_{\i}}^\nu$ is always irreducible as was shown in \secref{rk1} by an
explicit calculation.

Now let $\la\in\Lambda^{++}$ and $\nu\in\Lambda$ be such that $\nu\neq
w_0(\la)$, i.e., so that $\Gr^\la \cap S^\nu$ has positive
dimension. Let $K$ be an irreducible component of dimension $\langle
\la+\mu,\rhoc\rangle$ of $\Gr^\la \cap S^\nu$.

The following result is proved in \cite{BG}:

\sssec{\bf Proposition.}  \label{dominant} {\em For $\la,\nu$ and $K$
as above, there exist $\i\in\I$ and $\mu\in\Lambda_{M_\i}^{++}$ such
that an open dense subset of $K$ projects dominantly onto
$\Gr_{M_{\i}}^\mu \cap S_{M_{\i}}^\nu$ under the map $p_{\i}^\nu$ and
$\nu \neq w_{0,M_\i}(\mu)$.}

\smallskip

\noindent(In the above formula $w_{0,M_\i}$ is the ``longest'' element of
the Weyl group of $M_\i$, i.e. the $\i$-th simple reflection.)

\sssec{}

Recall from Sects.~\ref{groupcharacters} and \ref{admcharacters} that
the character $\chi_{-\nu}:N(\K_x)\to \GG_a$ is by definition a sum of
characters $^\i\chi_{\nu}$, $\i\in\I$. We can define functions
$^\i\chi^{\nu}_{-\nu}: S^{\nu}\to\GG_a$ so that
$$\chi^{\nu}_{-\nu}=\sum_{\i\in\I} \, {}^\i\chi^{\nu}_{-\nu}.$$

Let $\chi_{-\nu,M_\i}$ and $\chi^{\nu}_{-\nu,M_\i}$ be the corresponding
objects for $M_\i$.

\sssec{\bf Lemma.}
$^\i\chi^{\nu}_{-\nu}=\chi^{\nu}_{-\nu,M_\i} \circ p_{\i}^\nu$.

\bigskip

Combining this lemma with \propref{dominant} and the fact that
\propref{components} is true for $M_\i$ (proved in \secref{rk1}), we
obtain:

\sssec{\bf Corollary.} \label{nontrivrkone}
{\em If $\nu \neq w_0(\la)$, then for each component $K$ of $\Gr^\la
\cap S^{\nu}$, there exists $\i \in \I$, such that the restriction of
$^{\i}\chi^{\nu}_{-\nu}$ to $K(\Fqb)$ is non--constant.}

\sssec{}  \label{Taction}

The subgroup $T\subset T(\O_x)$ acts on $\Gr$. This action preserves
both $S^{\nu}$ and $\Gr^\la$. Since $T$ is connected, $K$ is preserved
as well.

The group $T$ also acts on $N(\K_x)$ by conjugation, and therefore on
the set of $N(\K_x)$--characters: for $\tau\in T$, $\chi \to \chi^\tau$,
where $\chi^\tau(n) = \chi(\tau \cdot n \cdot \tau^{-1})$. We have:
$$(^{{\mathfrak j}}\chi_{-\nu})^\tau = \alphach_{{\mathfrak j}}(\tau)
\cdot {}^{{\mathfrak j}}\chi_{-\nu}.$$

Without loss of generality we can assume that for all ${\mathfrak
j}\in\I$ the function $^{{\mathfrak j}}\chi_{-\nu}^{\nu}$ satisfies
$^{{\mathfrak j}}\chi_{-\nu}^{\nu}(\nu(t)) = 0$. Then, since the point
$\nu(t)\in\Gr$ is $T$--stable, we obtain:
\begin{equation}  \label{actonfunct}
(^{{\mathfrak j}}\chi^{\nu}_{-\nu})^\tau = \alphach_\jj(\tau) \cdot
{}^{{\mathfrak j}}\chi^{\nu}_{-\nu}.
\end{equation}

\bigskip

Now we can prove \propref{components}. It suffices to show that on the
level of $\Fqb$--points, the restriction of $\chi^{\nu}_{-\nu}$ to a
component $K$ as above is non--constant.  Suppose, to the contrary,  that
it is constant.

Letting $\i$ be as in \corref{nontrivrkone}, we obtain:
\begin{equation}    \label{zero}
^{\i}\chi_{-\nu}^{\nu}|_K= - \sum_{\jj \neq \i}{} \;\;
^{\jj}\chi_{-\nu}^{\nu}|_K+\on{const}
\end{equation}

Let $\tau$ be an element of $T$ satisfying $\alphach_{\i}(\tau)\neq
1$, $\alphach_{{\mathfrak j}}(\tau)=1$ for ${\mathfrak j} \neq \i$.
Apply $\tau$ to both sides of \eqref{zero}. According to
\eqref{actonfunct} we obtain:
$$\alphach_\i(\tau) \cdot {}^{\i}\chi_{-\nu}^{\nu}|_K =
\sum_{{\mathfrak j} \neq \i}{} \;\;
{}^{\jj}\chi_{-\nu}^{\nu}|_K+\on{const},$$
which, together with \eqref{zero}, contradicts the fact that
${}^{\i}\chi_{-\nu}^{\nu}|_K$ is not constant.

This completes the proof of \propref{components}.

\section{Proof of the conjecture from \cite{previous}}  \label{conj}

\ssec{The statement of the conjecture and its generalization}

\sssec{}

As was explained in the Introduction, one of the main motivations for
this paper was the following:

\sssec{\bf Conjecture.}(\cite{previous})   \label{spec}
{\it The cohomology
$H_c^k(\ol{\Gr}^\la\cap S^\nu, \A_\la|_{\ol{\Gr}^\la\cap S^\nu}\otimes
\chi^\nu_0|_{\ol{\Gr}^\la\cap S^\nu}{}^*(\J_\psi))$
vanishes unless $k=\langle2\nu,\rhoc\rangle$ and $\nu=\la$. In the
latter case it is canonically isomorphic to $\Ql$.}

\bigskip

In this section we will prove \thmmref{gencs1}, which is a
generalization of this conjecture. Recall that Theorem 1 states that for
$\la\in\Lambda^{++}$ and $\mu,\nu\in\Lambda$ with
$\mu+\nu\in\Lambda^{++}$ the cohomology
\begin{equation} \label{gencs2}
H_c^k(\ol{\Gr}^\la\cap S^{\nu}, \A_\la|_{\ol{\Gr}^\la\cap S^\nu}\otimes
\chi^\nu_{\mu}|_{\ol{\Gr}^\la\cap S^\nu}{}^*(\J_\psi))
\end{equation}
vanishes unless $k=\langle 2\nu,\rhoc\rangle$ and
$\mu\in\Lambda^{++}$. In the latter case, this cohomology identifies
canonically with $\Hom_{\GL}(V^{\la}\otimes V^{\mu},V^{\mu+\nu})$.

\sssec{}

As explained below, this theorem follows from
\corref{decomp} and \thmmref{one}(2) as a combination of two
facts. First, we can readily recognize the above cohomology groups as
the stalks of the LHS of \eqref{decomp1} (up to replacing $\la$ by
$-w_0(\la)$). Second, we can compute the stalks of the RHS
of \eqref{decomp1} using \thmref{one}.

\sssec{}

Note that \conjref{spec} is a special case of \thmmref{gencs1},
when we set $\mu=0$.

Now let $\mu$ be very large compared to $\la$ and $\nu$. In this case
the function $\chi^\nu_{\mu}:\ol{\Gr}^\la\cap S^\nu\to\GG_a$
is constant, i.e., $\chi^\nu_{\mu}|_{\ol{\Gr}^\la\cap
S^\nu}{}^*(\J_\psi)$ is isomorphic to $\Ql$. Note also that in the
case when $\la,\nu\ll\mu$, the vector space $\Hom_{\GL}(V^{\mu}\otimes
V^\la,V^{\nu+\mu})$ can be naturally identified with the dual of the
$\nu$--weight space $V^\la(\nu)$. Thus, \thmmref{gencs1} yields
the following result which was previously proved by Mirkovi\'c and one
of the authors by different methods:

\medskip

\noindent{\bf Theorem.}
{\em The cohomology
$H_c^k(\ol{\Gr}^\la\cap S^{\nu}, \A_\la|_{\ol{\Gr}^\la\cap S^\nu})$
vanishes unless $k=\langle 2\nu,\rhoc\rangle$ and, in the latter case,
it is isomorphic to $V^\la(\nu)^*$.}

\bigskip

Finally, note that \corref{corcomp} is also a special case of
\thmmref{gencs1} with $\mu=-\nu$ (see also Remark 2 below).

\sssec{} \label{remthmone}

Let us make several remarks concerning the structure of
\thmmref{gencs1}.

\noindent{\em Remark 1.}

The vanishing part of \thmmref{gencs1} is obvious when $\mu$ is
non--dominant. Indeed, the group $N(\O_x)$ acts on $\ol{\Gr}^\la\cap
S^{\nu}$ and the complex $\A_\la|_{\ol{\Gr}^\la\cap S^\nu}\otimes
\chi^\nu_{\mu}|_{\ol{\Gr}^\la\cap S^\nu}{}^*(\J_\psi)$ is
$(N(\O_x),\chi_{\mu})$--equivariant, while
$\chi_{\mu}|_{N(\O_x)}$ is non--trivial if $\mu$ is
non--dominant.

\medskip

\noindent{\it Remark 2.}

It follows immediately from \propref{dimestimate} that
$$H_c^k(\ol{\Gr}^\la\cap S^{\nu}, \A_\la|_{\ol{\Gr}^\la\cap
S^\nu}\otimes \chi^\nu_{\mu}|_{\ol{\Gr}^\la\cap
S^\nu}{}^*(\J_\psi))=0$$ if $k>\langle 2\nu,\rhoc\rangle$ and that for
$k=\langle 2\nu,\rhoc\rangle$ it is isomorphic to
$$H_c^{\langle \la+\nu,2\rhoc\rangle} (\Gr^\la\cap
S^\nu,\chi^\nu_{\mu}|_{\Gr^\la\cap S^\nu}{}^*(\J_\psi)).$$

\medskip

\noindent{\it Remark 3.}

Assume now that $G$ is defined and split over $\Fq$. Then the 
cohomology $$H_c^{\langle \nu,2\rhoc\rangle}(\ol{\Gr}^\la\cap S^{\nu},
\A_\la|_{\ol{\Gr}^\la\cap S^\nu}\otimes
\chi^\nu_{\mu}|_{\ol{\Gr}^\la\cap S^\nu}{}^*(\J_\psi))$$ is
also defined over $\Fq$. The action of the Frobenius on the
cohomology is easy to recover:

The previous remark implies that $\on{Fr}$ acts on it by $q^{\langle
\nu,2\rhoc\rangle}$, i.e., as on $\Ql(-\langle \nu,2\rhoc\rangle)$.

\medskip

\noindent{\it Remark 4.}

Consider the vector space $\Hom_{\GL}(V^\la\otimes V^{-w_0(\mu+\nu)},
V^{-w_0(\mu)})$. The previous
discussion implies that it can be naturally identified with
$$\Hom_{\GL}(V^\la\otimes V^\mu,V^{\nu+\mu})\simeq H_c^{\langle
\la+\nu,2\rhoc\rangle} (\Gr^\la\cap
S^\nu,\chi^\nu_{\mu}|_{\Gr^\la\cap S^\nu}{}^*(\J_\psi)).$$

Therefore, it acquires a basis labeled by the set $\{K\}$ of those
irreducible components $K$ of $\Gr^\la\cap S^\nu$, for which
$\chi^\nu_{\mu}|_K$ is non--constant (see the proof of
\corref{corcomp}).

It follows from \thmmref{sheaf-Satake} that
$\Hom_{\GL}(V^\la\otimes V^{-w_0(\nu+\mu)},
V^{-w_0(\mu)})$ admits another basis of geometric origin as follows.
Recall the scheme $\Conv$ (cf. \secref{intrconv}). For $\la_1$ and
$\la_2$ denote by $\Conv^{\la_1,\la_2}$ the corresponding locally
closed subscheme in $\Conv$. It is a fibration over $\Gr^{\la_1}$ with
a typical fiber isomorphic to $\Gr^{\la_2}$. It is known (\cite{MV}) that
the projection $p_2$ is a semi--small map $\Conv^{\la_1,\la_2}\to
\Gr^{\la_1+\la_2}$. Thus, to a triple of coweights $\la_1,\la_2,\la_3$
one can attach the set $K'(\la_1,\la_2,\la_3)$ of irreducible
components of the fiber of $\Conv^{\la_1,\la_2}$ over a point of
$\Gr^{\la_3}$. Hence the set
$K'(\la_1,\la_2,\la_3)$ parametrizes a basis for
$\Hom_{\GL}(V^{\la_1}\otimes V^{\la_2}, V^{\la_3})$.

Put $\la_1=\la$, $\la_2=-w_0(\mu+\nu)$, $\la_3=-w_0(\mu)$. One can
show that the sets $\{ K'(\la_1,\la_2,\la_3) \}$ and $\{K\}$ are in a
natural bijection and that the corresponding bases of the vector
space $\Hom_{\GL}(V^\la\otimes V^{-w_0(\mu+\nu)}, V^{-w_0(\mu)})$
coincide.

\ssec{Proof of \thmmref{gencs1}}

\sssec{}

First, let us show that it is enough to consider the case when
$[G,G]$ is simply--connected.

Let $G_1\twoheadrightarrow G$ be a surjection of reductive groups such
that its kernel is a (connected) torus. Then we have a surjection of
the corresponding lattices $\Lambda_1\twoheadrightarrow\Lambda$.  Fix
$\la_1,\mu_1,\nu_1\in\Lambda_1$ as in \thmmref{gencs1} and let
$\la,\mu,\nu$ be their images in $\Lambda$.

Starting with an arbitrary $G$, we can always find a group $G_1$
as above such that $[G_1,G_1]$ is simply--connected.
\thmmref{gencs1} for $G$ follows from \thmmref{gencs1} for
$G_1$ because of the following obvious fact:

\sssec{\bf Lemma.}  {\em The cohomology appearing in
\thmmref{gencs1} for $G_1$ and $\la_1,\mu_1,\nu_1$ is naturally
isomorphic to that of $G$ and $\la,\mu,\nu$.}

\sssec{}

Now we assume that $[G,G]$ is simply--connected. Fix a data of
$\varpi$ satisfying $\on{cond}_x(\varpi)=0$. As was explained before,
\thmmref{gencs1} will be obtained by computing explicitly the
restriction of $\Psib^{x,\mu+\nu}_\varpi\star\A_{-w_0(\la)}$ to
$_{x,\mu}\BunNft$. According to Remark 1 above, we can assume that $\mu$ is
dominant.

On the one hand, we claim that we have a canonical isomorphism
\begin{equation} \label{firstcomp}
j_\mu^*(\Psib^{x,\mu+\nu}_\varpi\star\A_{-w_0(\la)})\simeq
\Psi^{x,\mu}_\varpi\otimes
\Hom_{\GL}(V^\la\otimes V^\mu,V^{\nu+\mu}).
\end{equation}

Indeed, by \corref{decomp},
$$\Psib^{x,\mu+\nu}_\varpi\star\A_{-w_0(\la)}
\simeq \underset{\mu'}\oplus \Psib^{x,\mu'}_\varpi \otimes
\Hom_{\GL}(V^{\mu'},V^{-w_0(\la)}\otimes V^{\nu+\mu}).$$

But by \thmmref{one}(1), $j_\mu^*(\Psib^{x,\mu'}_\varpi)=0$ unless
$\mu'=\mu$. Therefore,
$$j_\mu^*(\Psib^{x,\mu+\nu}_\varpi\star\A_{-w_0(\la)})\simeq
\Psi^{x,\mu}_\varpi\otimes \Hom(V^{\mu},V^{-w_0(\la)}\otimes
V^{\nu+\mu})\simeq \Psi^{x,\mu}_\varpi\otimes \Hom(V^\la\otimes
V^\mu,V^{\nu+\mu}).$$

\sssec{}

On the other hand, let us compute
$j_\mu^*(\Psib^{x,\mu+\nu}_\varpi\star\A_{-w_0(\la)})$ using the
definition of the Hecke functors.

Recall the notation of \secref{stratifications}.  As
$\Psib^{x,\mu+\nu}_\varpi\simeq \Psi^{x,\mu+\nu}_\varpi{}_{!}$ (here
we use \thmmref{one}(1) once again), using base change,
$j_\mu^*(\Psib^{x,\mu+\nu}_\varpi\star\A_{-w_0(\la)})$ can be computed
as follows:

We  $*$--restrict
$\Psib_\varpi^{x,\mu+\nu}\tboxtimes\A_{-w_0(\la)}$ to
$Z^{\mu,\mu+\nu}$ and then apply the $!\,$-push--forward with respect to the
map $'\hl:Z^{\mu,\mu+\nu}\to{}_{x,\mu}\BunNft$.
Using \lemref{partdescr}(1), we identify $Z^{\mu,\mu+\nu}$ with the
fibration
$$_\mu\widetilde{\mathcal N}_x^{\epsilon_\mu} \underset{N(\O_x)}\times
S^{\nu}\overset{'\hl}\longrightarrow{}_{x,\mu}\BunNftw$$ and the
support of
$(\Psib_\varpi^{x,\mu+\nu}\tboxtimes\A_{-w_0(\la)})|_{Z^{\mu,\mu+\nu}}$
is the substack $$_\mu\widetilde{\mathcal
N}_x^{\epsilon_\mu}\underset{N(\O_x)}\times (\ol{\Gr}^\la\cap
S^{\nu})$$ of $Z^{\mu,\mu+\nu}$.

Moreover, according to formula \eqref{left-right} and
\lemref{functionsmatch}(2), the restriction of
$$\Psib_\varpi^{x,\mu+\nu}\tboxtimes\A_{-w_0(\la)}$$ to this substack
can be identified with the sheaf
$$(\Psi_\varpi^{x,\mu})\tboxtimes(\A_\la|_{\ol{\Gr}^\la\cap S^\nu} \otimes
\chi^\nu_{\mu}|_{\ol{\Gr}^\la\cap S^\nu}{}^*(\J_\psi)) 
[\langle\nu,2\rhoc\rangle]$$ (note that
$\langle\nu,2\rhoc\rangle=\on{dim}({}_{x,\mu+\nu}\BunNft) -
\on{dim}({}_{x,\mu}\BunNft)$ and that is why it enters into the above
formula.)

Therefore, by the projection formula, we obtain that
$j_\mu^*(\Psib^{x,\mu+\nu}_\varpi\star\A_{-w_0(\la)})$ is isomorphic
to
$$\Psi_\varpi^{x,\mu}\otimes H_c^{\bullet}(\A_\la|_{\ol{\Gr}^\la\cap
S^\nu}\otimes \chi^\nu_{\mu}|_{\ol{\Gr}^\la\cap
S^\nu}{}^*(\J_\psi)[\langle\nu,2\rhoc\rangle]).$$ Comparing it with
\eqref{firstcomp}, we obtain the statement of \thmmref{gencs1}.

\end{document}